\documentclass[a4paper, 12pt]{article}
\usepackage{latexsym}
\usepackage{amsmath}
\usepackage{amssymb}
\usepackage{amsthm}
\usepackage{amscd}
\usepackage[dvipdfmx]{graphicx}

\renewcommand{\thefootnote}{}

\begin{document}
\title{ {\bf On a certain nilpotent extension  over $\mathbb{Q}$ of degree 64 and the 4-th multiple residue symbol}}
\author{Fumiya Amano }
\date{ }
\maketitle

\footnote{{\footnotesize 2010 Mathematics Subject Classification Primary 11A15 11R32, Secondary 57M27. \\
$\;\;\;\;$ Key words: R\'{e}dei symbol, Milnor invariant, 4-th multiple residue symbol}}
\renewcommand{\thefootnote}{*1}
{\small {\bf Abstract}. In this paper, we introduce the 4-th multiple residue symbol $[p_1,p_2,p_3,p_4]$ for certain four prime numbers $p_i$'s, which extends the Legendre symbol $\left(\frac{p_1}{p_2}\right)$ and the R\'{e}dei triple symbol $[p_1,p_2,p_3]$ in a natural manner.
For this we construct concretely a certain nilpotent extension $K$ over $\mathbb{Q}$ of degree $64$, where ramified prime numbers are 
$p_1, p_2$ and $p_3$, such that the symbol $[p_1,p_2,p_3,p_4]$ describes the decomposition law of $p_4$ in the extension $K/\mathbb{Q}$.
We then establish the relation of our symbol $[p_1,p_2,p_3,p_4]$ and the 4-th arithmetic Milnor invariant $\mu_2(1234)$ (an arithmetic analogue of the 4-th order linking number) by showing $[p_1,p_2,p_3,p_4] = (-1)^{\mu_2(1234)}$.}

\quad \\
\begin{center}
\bf{Introduction } 
\end{center}

\quad As is well known, for two odd prime numbers $p_1$ and $p_2$, the Legendre symbol $\left( \frac{p_1}{p_2} \right) $ describes the decomposition law of $p_2$ in the quadratic extension $\mathbb{Q}(\sqrt{p_1})/\mathbb{Q}$.
In 1939, L. R\'{e}dei ([R]) introduced a triple symbol with the intention of a generalization of the Legendre symbol and Gauss' genus theory.
For three prime numbers  $p_i \equiv 1 \pmod{4} \ (i = 1, 2, 3)$  with $\left( \frac{p_i}{p_j} \right) = 1 \ (1 \le i \neq j \le 3)$, the R\'{e}dei triple symbol [$p_1, p_2, p_3 $] describes the decomposition law of $p_3$ in a Galois extension over $\mathbb{Q}$ where all ramified prime numbers are $p_1$ and $p_2$ and the Galois group is the dihedral group $D_8$ of order 8.

Although a meaning of the R\'{e}dei symbol had been obscure for a long time, in 2000, M. Morishita ([Mo1,2,3]) interpreted the R\'{e}dei symbol as an arithmetic analogue of a mod 2 triple linking number, following the analogies between knots  and primes.
In fact, he introduced arithmetic analogue $\mu_{2}(12 \cdots n) \in \mathbb{Z}/2\mathbb{Z}$ of Milnor's link invariants (higher order linking numbers) 
for prime numbers $p_1, \cdots, p_n$ such that
\[
\left( \frac{p_1}{p_2} \right) = (-1)^{\mu_{2}(12)}, \ [p_1, p_2, p_3] = (-1)^{\mu_{2}(123)}.
\]

Since it is difficult to compute arithmetic Milnor invariants by the definition, it is desirable to construct Galois extensions $K_n/\mathbb{Q}$ concretely such that $[p_1, \ldots, p_n] = (-1)^{\mu_{2}(12 \cdots n)}$ describes the decomposition law of $p_n$ in $K_n/\mathbb{Q}$, just as in the cases of the Legendre symbol where $K_2$ is a quadratic extension and the R\'{e}dei triple symbol where $K_3$ is a dihedral extension of degree 8.
As we shall explain in Section 2.1, link theory suggests that the desired extension $K_n/\mathbb{Q}$ should be a Galois extension such that all ramified prime numbers are $p_1, \ldots, p_{n-1}$ and the Galois group is the nilpotent group
\[
N_{n}(\mathbb{F}_2) = \left\{ \left.  \begin{aligned}  \left (
\begin{array}{cccc}
1 & *  & \cdots & *   \\ 
0 & 1 &  \ddots  & \vdots \\
\vdots & \ddots & \ddots & * \\
0 & \cdots & 0 & 1 \\
\end{array}
\right) \end{aligned} \right| * \in \mathbb{F}_2 \right\}
\]
consisting of $n \times n$ unipotent upper-triangular matrices over $\mathbb{F}_2$.
Note that 
$N_{2}(\mathbb{F}_2) = \mathbb{Z}/2\mathbb{Z}$ and $N_{3}(\mathbb{F}_2) = D_8$.

The purpose of this paper is to construct concretely such an extension $K_n/\mathbb{Q}$ for $n = 4$ in a natural manner extending R\'{e}dei's dihedral extension.
We then introduce the 4-th multiple residue symbol $[p_1, p_2, p_3, p_4]$ 
describing the decomposition law of $p_4$ in $K_4/\mathbb{Q}$ and prove that it coincides 
with the 4-th Milnor invariant $\mu_{2}(1234)$,
\[
[p_1, p_2, p_3, p_4] = (-1)^{\mu_{2}(1234)}.
\]
\\ 
\quad {\it Notation}. For a number field $k$, we denote by ${\cal O}_k$ the ring of integers of $k$.
For a group $G$ and $d \in \mathbb{N}$, we denote by $G^{(d)}$ 
the $d$-th term of the lower central series of $G$ defined by
$G^{(1)} := G, G^{(d+1)} := [G,G^{(d)}]$.
For a ring $R$, $R^{\times}$ denotes the group of invertible elements
of $R$. 
\\

\begin{center}
{\bf \S 1\quad R\'{e}dei's dihedral extension and triple symbol}
\end{center}

\quad In this section, we recall the construction of R\'{e}dei's dihedral extension 
and triple symbol ([R]), which will be used later.
We also give some basic 
properties of R\'{e}dei's dihedral extension and triple symbol.
\\

{\bf 1.1. The R\'{e}dei extension.}
\quad Let $p_1$ and $p_2$ be distinct prime numbers satisfying
\[
p_i \equiv 1 \pmod{4} \ (i = 1, 2), \quad \left( \frac{p_1}{p_2} \right) = \left( \frac{p_2}{p_1} \right) = 1  \leqno{(1.1.1)}
\]
We set $k_i =  \mathbb{Q}(\sqrt{p_i}) \ (i = 1, 2)$.
It follows from this assumption (1.1.1) that we have the following Lemma.
\\
\\
\quad {\bf Lemma 1.1.2} ([A, Lemma 1.1]){\bf .}
 {\it There are integers $x, y, z$ satisfying the following conditions$:$
\\
$(1) \;\; x^2 - p_1y^2 - p_2z^2  = 0.	$
\\
$(2) \;\; \mathrm{g.c.d}(x,y,z) = 1,\quad  y \equiv 0 \pmod{2},\quad x - y \equiv 1 \pmod{4}.$ \\
Furthermore, for a given prime ideal $\frak{p}_2$ of ${\cal O}_{k_1}$ lying over $p_2$, we can find integers $x, y, z$ which satisfy $(1), (2)$ and $(x + y\sqrt{p_1}) = \frak{p}_2^m$ for an odd positive integer $m$.}
\\
\\
Let $\boldsymbol{a} = (x,y,z)$ be a triple of integers satisfying the conditions (1), (2) in Lemma 1.1.2.
Then let $\alpha = x + y\sqrt{p_1}$ and set
\[
k_{\boldsymbol{a}} = \mathbb{Q}(\sqrt{p_1},\sqrt{p_2},\sqrt{\alpha}).
\leqno{(1.1.3)}
\]
The following theorem was proved by L. R\'{e}dei ([R]).
\\
\\
\quad {\bf Theorem 1.1.4} ([R]){\bf .}
 {\it 
$(1)$ The field $k_{\boldsymbol{a}}$ is a Galois extension over $\mathbb{Q}$ whose Galois group is the dihedral group of order $8$.
\\
$(2)$ Let $d(k_1(\sqrt{\alpha})/k_1)$ be the relative discriminant of the extension  $k_1(\sqrt{\alpha})/k_1$.
Then we have $N_{k_1/\mathbb{Q}}(d(k_1(\sqrt{\alpha})/k_1)) = (p_2)$.
In particular, all prime numbers ramified in $k_{\boldsymbol{a}}/\mathbb{Q}$ are $p_1$ and $p_2$ with ramification index $2$.
}
\\

The fact that $k_{\boldsymbol{a}}$ is independent of the choice of ${\boldsymbol{a}}$ was also shown in [R].
The author gave an alternative proof of this fact in [A], based on a proof communicated by D. Vogel ([V2]).
\\
\\
\quad {\bf Theorem 1.1.5} ([A, Corollary 1.5]){\bf .}  
{\it 
 A field $k_{\boldsymbol{a}}$ is independent of the choice of $\boldsymbol{a} = (x, y, z)$ satisfying $(1)$ and $(2)$ in Lemma $1.1.2$, namely, it depends only on a set $ \{ p_1, p_2 \} $.
}
\\
\\
\quad {\bf Definition 1.1.6.}  
By Proposition 1.1.5, we denote by $k_{ \{ p_1, p_2 \} }$ the field $k_{\boldsymbol{a}} = \mathbb{Q}(\sqrt{p_1},\sqrt{p_2},\sqrt{\alpha})$ given by (1.1.3) and call  $k_{ \{ p_1, p_2 \} }$ the {\it R\'{e}dei extension} over $\mathbb{Q}$ associated to a set $\{ p_1, p_2\} $ satisfying (1.1.1).
\\
\\
The following theorem shows that the R\'{e}dei extension $k_{ \{ p_1, p_2 \} }/\mathbb{Q}$ is characterized by the information on the Galois group and the ramification given in Theorem 1.1.4.
\\
\\
\quad {\bf Theorem 1.1.7} ([A, Theorem 2.1]){\bf .}
{\it 
Let $p_1$ and $p_2$ be prime numbers satisfying the condition $(1.1.1)$.
Then the following conditions on a number field $K$ are equivalent$:$
\\
$(1)$ $K$ is the R\'{e}dei extension $k_{ \{ p_1, p_2 \} }$.
\\
$(2)$ $K$ is a Galois extension over $\mathbb{Q}$ such that the Galois group is the dihedral group $D_8$ of order $8$ and prime numbers ramified in $K/\mathbb{Q}$ are $p_1$ and $p_2$ with ramification index $2$.
}
\\

\quad {\bf 1.2. The R\'{e}dei triple symbol.}  Let $p_1, p_2$ and $p_3$ be three prime numbers satisfying 
\[
p_i \equiv 1 \pmod{4} \ (i = 1, 2, 3), \ \left( \frac{p_i}{p_j} \right) = 1 \ (1 \le i \neq j \le 3). \leqno{(1.2.1)}
\] 
Let $k_{ \{ p_1, p_2\} }$ be the R\'{e}dei extension over $\mathbb{Q}$ associated to a set $\{ p_1, p_2\}$ (Definition 1.1.6).
\\
\\
\quad {\bf Definition 1.2.2.}  
We define {\it R\'{e}dei triple symbol} $[p_1, p_2, p_3]$ by
\[
[p_1, p_2, p_3] = \left\{
\begin{array}{rl}
1 &\quad \mbox{if  $p_3$ is completely decomposed in}\  k_{ \{ p_1, p_2\} }/\mathbb{Q} ,\\
-1 &\quad \mbox{otherwise}.
\end{array}
\right.
\]
\\
The following theorem is a reciprocity law for the R\'{e}dei triple symbol:
\\
\\
\quad {\bf Theorem 1.2.3} ([R], [A, Theorem 3.2]){\bf .}  {\it  We have }
\[
[p_1, p_2, p_3] = [p_i, p_j, p_k].
\]
{\it for any permutation $\{ i, j, k \} $ of $\{ 1, 2, 3 \}$.}
\\

\begin{center}
{\bf \S 2 \quad Milnor invariants}
\end{center}

 In this section, we recall the arithmetic analogues of Milnor invariants 
of a link introduced by M. Morishita ([Mo1,2,3]) and clarify a meaning of the R\'{e}dei extension and the R\'{e}dei triple symbol in Section 1 from the viewpoint of the analogy  between knot theory and number theory.
The underlying idea is based on the following analogies between knots and primes (cf. [Mo4]):

\begin{center}
	\begin{tabular}{ | c | c |} \hline
		knot & prime \\
		${\cal K} : S^1 \hookrightarrow \mathbb{R}^3$ & Spec$(\mathbb{F}_p) \hookrightarrow \mbox{Spec}(\mathbb{Z})$ \\   \hline
		link & finite set of primes \\
		${\cal L} = {\cal K}_1 \cup \cdots  \cup {\cal K}_r $ & $S = \{ p_1, \ldots, p_r \}$ \\   \hline
		$X_{ {\cal L} } = \mathbb{R}^3 \setminus {\cal L} $ & $X_S = \mbox{Spec}(\mathbb{Z})\setminus S$	\\   \hline	
		   & Galois group with restricted ramification \\
		link group & $G_S = \pi_{1}^{ \mbox{\footnotesize{\'{e}t}} } (X_S ) = \mbox{Gal}(\mathbb{Q}_S/\mathbb{Q})$ \\
		$G_{\cal L} = \pi_1(X_{ {\cal L} })$  & $\mathbb{Q}_S$ : maximal extension over $\mathbb{Q}$	\\
		   &  unramified outside $S \cup \{ \infty \}$   \\   \hline
	\end{tabular}
\end{center}

In the following, we firstly explain Milnor invariants of a link and their 
meaning in nilpotent coverings of $S^3$ ([Mi2], [Mu]). We then discuss their arithmetic analogues for prime numbers where the R\'{e}dei triple symbol is interpreted as an arithmetic analogues of a triple Milnor invariant. The analogy also suggests that a natural generalization of the Legendre and R\'{e}dei symbols, called a multiple residue symbol $[p_1,\dots ,p_n]$, should describe the decomposition law of $p_n$ in a certain nilpotent extension over $\mathbb{Q}$ unramified outside $p_1,\dots, p_{n-1}$ and $\infty$ ($\infty$ being the infinite prime).\\

{\bf 2.1. Milnor invariants of a link.}   Let ${\cal L} = {\cal K}_1 \cup \cdots  \cup {\cal K}_r$ be a link with $r$ components in $\mathbb{R}^3$ and let $X_{ {\cal L} } = \mathbb{R}^3 \setminus {\cal L}$ and $G_{ {\cal L} } := \pi_1(X_{ {\cal L} })$ be the link group of ${\cal L}$.
Let $F$ be the free group on the words $x_1 , \ldots , x_r$ where $x_i$ represents a meridian of ${\cal K}_i$.
The following theorem is due to J. Milnor.
\\
\\
\quad {\bf Theorem 2.1.1} ([Mi2, Theorem 4]){\bf .}
{\it \quad
For each $d \in \mathbb{N}$, there is $y_{i}^{(d)} \in F $ such that
\begin{align*}
G_{\cal L}/G_{ {\cal L} }^{(d)} &= \langle x_1, \ldots, x_r \mid [x_1,y_{1}^{(d)}] = \cdots = [x_r,y_{r}^{(d)}] = 1, \ F^{(d)} = 1 \rangle,  \\
y_{j}^{(d)} &\equiv y_{j}^{(d + 1)} \mod{F^{(d)}},
\end{align*}
where $y_{j}^{(d)}$ is a word representing  a longitude of ${\cal K}_j$ in $G_{\cal L}/G_{ {\cal L} }^{(d)}$.
}
\\

Let $\mathbb{Z} \langle\langle X_1, \ldots ,X_r\rangle\rangle$ be the algebra of non-commutative formal power series of variables $X_1, \ldots ,X_r$ over $\mathbb{Z}$, and let 
\[
M : F \longrightarrow \mathbb{Z}\langle\langle X_1, \ldots ,X_r\rangle\rangle^{\times}
\]
be the Magnus homomorphism defined by
\[
M(x_i) := 1 +X_i, \ M(x_{i}^{-1}) := 1 - X_i + X_{i}^2 - \cdots (1 \le i \le r). 
\]
For $f \in F$, $M(f)$ has the form
\[
M(f) = 1 + \sum_{n=1}^{\infty} \sum_{1 \le i_1, \ldots, i_n \le r} \mu (i_1 \cdots i_n ; f )  X_{i_1} \cdots X_{i_n},
\]
where the coefficients $\mu (i_1 \cdots i_n ; f )$ are called the {\it Magnus coefficients}.

Let $\mathbb{Z}[F]$ be the group algebra of $F$ over $\mathbb{Z}$ and let $\epsilon_{\mathbb{Z}[F]} : \mathbb{Z}[F] \rightarrow \mathbb{Z}$ be the augmentation map. We note that the Magnus coefficients can be written 
in terms of the Fox derivative introduced in [F]:
\[
\mu (i_1 \cdots i_n ; f ) = \epsilon_{\mathbb{Z}[F]} \left( \frac{\partial^{n} f }{\partial x_{i_1} \cdots \partial x_{i_n}} \right).
\]
\\

For the word $y_{j}^{(d)}$ in Theorem 2.1.1, we set 
\[
\mu^{(d)} (i_1 \cdots i_n j) := \mu (i_1 \cdots i_n ; y_{j}^{(d)}).
\]
Since $\mu (i_1 \cdots i_n ; f ) = 0$ for $f \in F^{(d)}$ if $d > n$, by Theorem 2.1.1,  $\mu^{(d)} (I)$ is independent of $d$ if $d \ge |I|$, where $|I|$ denotes the length of a multi-index $I$.
Define $\mu(I) := \mu^{(d)}(I) \; (d \gg 1)$.
For a multi-index $I$ with $|I| \ge 2$, we define $\Delta(I)$ to be the ideal of $\mathbb{Z}$ generated by $\mu(J)$ where $J$ runs over cyclic permutations of proper subsequences of $I$.
If $|I| = 1$, we set $\mu(I) : = 0$ and $\Delta(I) := 0$ .
The {\it Milnor $\overline{\mu}$-invariant} is then defined by
\[
\overline{\mu}(I) := \mu(I) \mod \Delta(I).
\]
The fundamental results, due to Milnor, are as follows.
\\
\\
\quad {\bf Theorem 2.1.2} ([Mi2, Theorems 5, 6]){\bf .}  
 {\it 
$(1)$ $\overline{\mu}(i j) = \mathrm{lk}({\cal K}_i, {\cal K}_j) \ (i \neq j)$. \\
$(2)$ If $2 \le |I| \le d$, $\overline{\mu}(I)$ is a link invariant of ${\cal L}$. \\
$(3)$ $(\mathrm{Shuffle}$ $\mathrm{relation} )$ For any $I, J (|I|, |J| \ge 1)$ and $ i \  (1 \le i \le r)$, we have 
\[
\sum_{H \in {\rm PSh}(I, J) } \overline{\mu}(Hi) \equiv 0 \mod  \mathrm{g.c.d} \{ \Delta(Hi) \mid H \in {\rm PSh}(I, J) \}
\]
where ${\rm PSh}(I, J)$ stands for the set of results of proper shuffles of $I$ and $J$ $(  \mathrm{cf. \ [CFL]} )$. \\
$(4)$ $(\mathrm{Cyclic}$ $\mathrm{symmetry} )$.
$\overline{\mu}(i_1 \cdots i_n) = \overline{\mu}(i_2 \cdots i_n i_1) = \cdots = \overline{\mu}(i_n i_1 \cdots i_{n - 1})$
}
\\
\\
\quad {\bf Example 2.1.3}.
For a multi-index $I$ $(|I| \ge 2)$, $\overline{\mu}(I) = \mu(I) $ is an integral link invariant if $\mu(J) =0$ for all multi-index $J$ with $ | J | < | I |$.
For example, let ${\cal L} = {\cal K}_1 \cup {\cal K}_2 \cup {\cal K}_3$ be the following {\it Borromean rings}:
\newpage 
\begin{center}
\includegraphics[height=4cm]{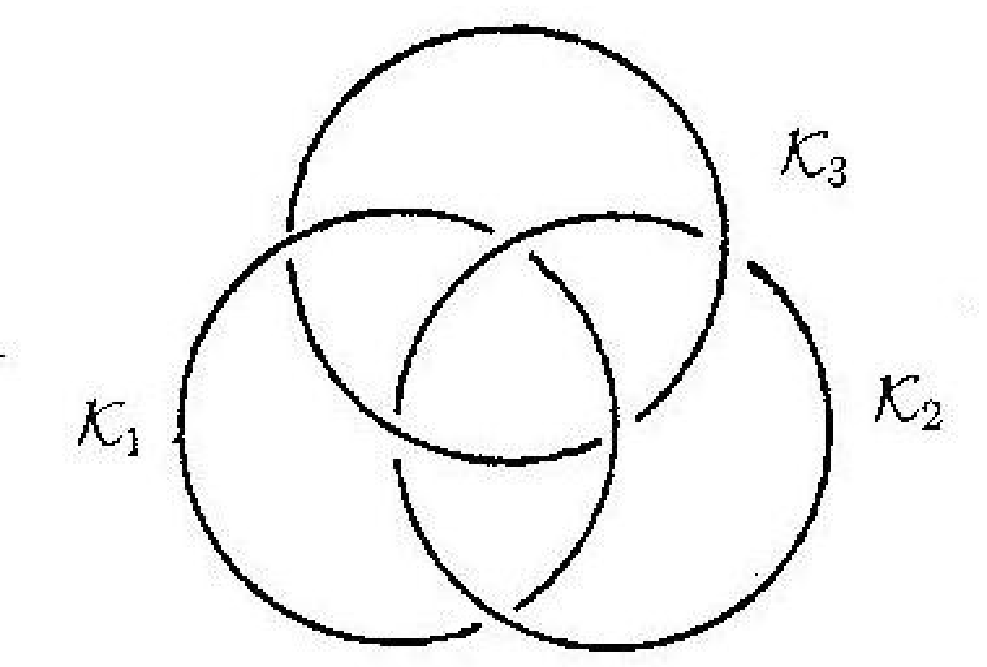}
\end{center}
Then $\mu(I) =0$ if $| I | \le 2$ and hence $\mu(I) \in \mathbb{Z}$ for $| I | = 3$.
In fact, we have $\mu(ijk) = \pm 1$ if $ijk$ is a permutation of 123 and $\mu(ijk) = 0$ otherwise.

Move generally, let ${\cal L} = {\cal K}_1 \cup \cdots  \cup {\cal K}_r$ be the following link, called the {\it Milnor link} ([Mi1, 5]).
\begin{center}
\includegraphics[height=4cm]{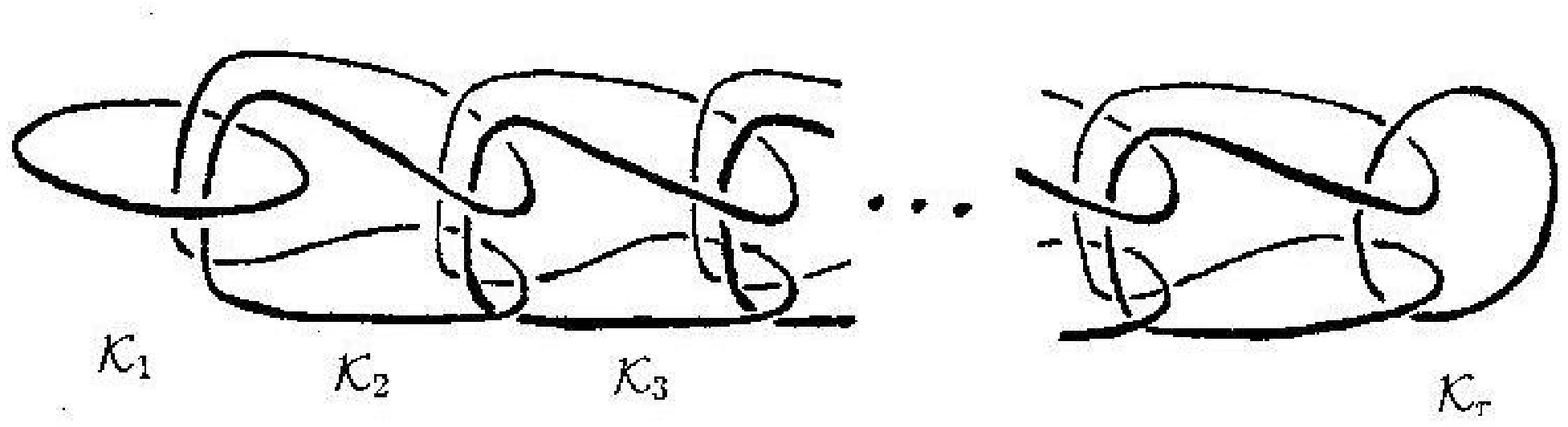}
\end{center}
We easily see that the link obtained by removing any one component ${\cal K}_i$ from ${\cal L}$ is trivial.
So $\mu(I) =0$ if $| I | \le {n - 1}$ and $\mu(I) \in \mathbb{Z}$ if $| I | = n$.
For instance, $\mu(12 \cdots n) = 1$.
\\

Next, we recall that Milnor invariants may be regarded as invariants associated to nilpotent coverings of $S^3$. For a commutative ring $R$, let $N_n(R)$ be the group consisting of $n$ by $n$ unipotent uppertriangular matrices.
For a multi-index $I = (i_1 \cdots i_n) (n \ge 2)$, we define the map $\rho_{I} : F \rightarrow N_{n}(\mathbb{Z}/\Delta(I))$ by
\[
\rho_{I}(f) := \left(
\begin{array}{ccccc}
	 1 & \epsilon(\frac{\partial f}{\partial x_{i_1}}) & \epsilon(\frac{\partial^{2} f}{\partial x_{i_1} \partial x_{i_2}}) & \cdots & \epsilon(\frac{\partial^{n - 1} f}{\partial x_{i_1} \cdots \partial{x_{i_{n-1}} } })\\
	 0 & 1 & \epsilon(\frac{\partial f}{\partial x_{i_2} }) & \cdots &\epsilon(\frac{\partial^{n - 2} f}{\partial x_{i_2} \cdots \partial{x_{i_{n-1}} } })\\
	 \vdots & \ddots & \ddots & \ddots & \vdots \\
	 \vdots & & \ddots & 1 & \epsilon(\frac{\partial f}{\partial x_{i_{n-1}} })\\
	 0 & \cdots & \cdots & 0 & 1 \\
\end{array}
\right)
\mod \Delta(I),
\]
where we set $\epsilon =  \epsilon_{\mathbb{Z}[F]}$ for simplicity. It can be shown by the property of the Fox derivative that $\rho_{I}$ is a homomorphism.
\\
\\
\quad {\bf Theorem 2.1.6} ([Mo4, Theorem 8.8], [Mu]){\bf .}  
 {\it 
$(1)$ The homomorphism $\rho_{I}$ factors through the link group $G_{\cal L}$.
Furthermore it is surjective if $i_1, \dots, i_{n - 1}$ are all distinct.  \\
$(2)$ Suppose that $i_1, \dots, i_{n - 1}$ are all distinct.
Let $X_I \rightarrow X_{\cal L}$ be the Galois covering corresponding to $\mathrm{Ker}(\rho_{I})$ whose Galois group $\mathrm{Gal}(X_I/X_{\cal L}) = N_{n}(\mathbb{Z}/\Delta(I))$.
When $\Delta(I) \neq 0$, let $M_I \rightarrow S^3$ be the Fox completion of $X_I \rightarrow X_{\cal L}$, a Galois covering ramified over the link ${\cal K}_{i_1} \cup \cdots  \cup {\cal K}_{i_{n -1}}$.
For a longitude $\beta_{i_n}$ of ${\cal K}_{i_n}$, one has

\[
\rho_{I}(\beta_{i_n}) = \left(
\begin{array}{ccccc}
	 1 & 0  & \cdots & 0 & \overline{\mu}(I)  \\ 
	 0 & 1 & \cdots &  & 0 \\
	 \vdots & \ddots & \ddots &  & \vdots \\
	 \vdots &  & \ddots & 1 & 0 \\
	 0 & \cdots & \cdots & 0 & 1 \\
\end{array}
\right)
\]
and hence the following holds$:$
\[
\overline{\mu}(I) = 0 \Longleftrightarrow {\cal K}_{i_n} is\ completely\ decomposed\ in\ M_I \rightarrow S^3.
\]
}

{\bf 2.2. Milnor invariants for prime numbers.}   Let $S = \{ p_1, \cdots, p_r \}$ be a 
set of $r$ distinct odd prime numbers and let $G_S := \pi_{1}^{ \footnotesize{ \mbox{\'{e}t} } }(\mathrm{Spec}(\mathbb{Z}) \setminus S)$.
In order to get the analogy of the link case, we consider the maximal pro-2 quotient, denoted by $G_S(2)$, of $G_S$ which is the Galois group of the maximal pro-2 extension $\mathbb{Q}_S(2)$ over $\mathbb{Q}$ which is unramified outside $S \cup \{ \infty \}$. Here we fix an algebraic closure $\overline{ \mathbb{Q} }$ of $\mathbb{Q}$ containing $\mathbb{Q}_S(2)$. We also fix an algebraic closure $\overline{ \mathbb{Q} }_{p_i}$ of $\mathbb{Q}_{p_i}$ and an embedding $\overline{ \mathbb{Q} } \hookrightarrow \overline{ \mathbb{Q} }_{p_i}$ for each $i$.
Let $\mathbb{Q}_{p_i}(2)$ be the maximal pro-$2$ extension of $\mathbb{Q}_{p_i}$ contained in $\overline{ \mathbb{Q} }_{p_i}$. Then we have
\[
\mathbb{Q}_{p_i}(2) = \mathbb{Q}_{p_i}( \zeta_{2^n}, \sqrt[2^n]{p_i} \mid n \ge 1)
\]
where $\zeta_{2^n} \in \overline{ \mathbb{Q} }$ is primitive $2^n$-th root of unity such that $\zeta_{2^t}^{2^s} = \zeta_{ 2^{t -s} }$ $(t \ge s)$.
The local Galois group Gal($\mathbb{Q}_{p_i}(2)/\mathbb{Q}_{p_i}$) is then topologically generated by the monodromy $\tau_i$ and the extension of the Frobenius automorphism $\sigma_i$ defined by 
\[
\left. \begin{array}{ll}
\tau_i (\zeta_{2^n}) = \zeta_{2^n}, &\quad \tau_i (\sqrt[2^n]{p_i}) =  \zeta_{2^n} \sqrt[2^n]{p_i}, \\
\sigma_i (\zeta_{2^n}) = \zeta_{2^n}^{p_1}, &\quad \sigma_i (\sqrt[2^n]{p_i}) = \sqrt[2^n]{p_i} 
\end{array}
\right. \leqno{(2.2.1)}
\]
and $\tau_i, \sigma_i$ are subject to the relation $\tau_i^{p_i -1}[\tau_i, \sigma_i] = 1$.

The embedding $\overline{ \mathbb{Q} } \hookrightarrow \overline{ \mathbb{Q} }_{p_i}$ induces the embedding $\mathbb{Q}_S(2) \hookrightarrow \mathbb{Q}_{p_i}(2)$ and hence the homomorphism $\eta_i : \mbox{Gal}(\mathbb{Q}_{p_i}(2)/\mathbb{Q}_{p_i})$ $ \rightarrow G_S$.
We denote by the same $\tau_i, \sigma_i$ the images of $\tau_i, \sigma_i$ under $\eta_i$.
Let $\hat{F}$ denote the free pro-2 group on the  words $x_1, \ldots ,x_r$ where $x_i$ represents $\tau_i$.
The following theorem, due to H. Koch, may be regarded as an arithmetic analogue of Milnor's Theorem 2.1.1.
\\
\\
\quad {\bf Theorem 2.2.2} ([K2, Theorem 6.2]){\bf .}
{\it \quad 
The pro-$2$ group $G_S(2)$ has the following presentation$:$
\[
G_S(2) = \langle x_1, \dots, x_r \mid x_{1}^{p_1 - 1}[x_1,y_{1}] = \cdots = x_{r}^{p_r - 1}[x_r,y_{r}] = 1 \rangle,
\]
where $y_{j} \in \hat{F}$ is the pro-$2$ word which represents $\sigma_j$.
}
\\

Set $e_S := \mathrm{max}\{ e \mid p_i \equiv 1 \mod{2^e} \ (1 \le i \le r ) \}$ and fix $m = 2^e$ ($1 \le e \le e_S$).
Let $\mathbb{Z}_{2}\langle\langle X_1, \ldots ,X_r\rangle\rangle$ be the algebra of non-commutative formal power series of variables $X_1, \ldots ,X_r$ over $\mathbb{Z}_2$, the ring of 2-adic integers, and let 
\\
\[
\hat{M} : \hat{F} \longrightarrow \mathbb{Z}_2\langle\langle X_1, \ldots ,X_r\rangle\rangle^{\times}
\]
be the pro-2 Magnus embedding ([K1, 4.2]).
For $f \in \hat{F}$, $\hat{M}(f)$ has the from
\[
\hat{M}(f) = 1 + \sum_{1 \le i_1, \ldots i_n \le r} \hat{\mu }(i_1 \cdots i_n ; f )  X_{i_1} \cdots X_{i_n},
\]
where the coefficients $\hat{\mu} (i_1 \cdots i_n ; f )$ are called the {\it $2$-adic Magnus coefficients}. We let 
$$
M_2 : \hat{F} \longrightarrow \mathbb{F}_2\langle\langle X_1, \ldots ,X_r\rangle\rangle^{\times}
$$
be the mod 2 Magnus embedding defined by composing $\hat{M}$ with the natural homomorphism $ \mathbb{Z}_2\langle\langle X_1, \ldots ,X_r\rangle\rangle^{\times} \longrightarrow \mathbb{F}_2\langle\langle X_1, \ldots ,X_r\rangle\rangle^{\times}$.

Let $\mathbb{Z}_2[[\hat{F}]]$ be the complete group algebra over $\mathbb{Z}_2$ and let $\epsilon_{\mathbb{Z}_2[[\hat{F}]]} : \mathbb{Z}_2[[\hat{F}]] \rightarrow \mathbb{Z}_2$ be the augmentation map.
In terms of the pro-2 Fox free derivative ([I], [O]), the 2-adic Magnus coefficients are written as
\[
\hat{\mu} (i_1 \cdots i_n ; f ) = \epsilon_{\mathbb{Z}_2[[\hat{F}]]} \left( \frac{\partial^{n} f }{\partial x_{i_1} \cdots \partial x_{i_n}} \right).
\]
\\
\quad For the word $y_j$ in Theorem 2.2.2, we set
\[
\hat{\mu} (i_1 \cdots i_n j) := \hat{\mu} (i_1 \cdots i_n ; y_j) 
\]
and we set, for a multi-index $I$,
\[
\mu_m(I) := \hat{\mu}(I) \mod m.
\]
For a multi-index with $I$ with $1\le |I| \le 2^{e_S}$, let $\Delta_m(I)$ be the ideal of $\mathbb{Z}/m\mathbb{Z}$ generated by $\binom {2^{e_S}}{t}$ $(1\le t \le |I|)$ and $\mu_m(J)$ ($J$ running over cyclic permutation of proper subsequences of $I$).
The Milnor $\overline{\mu}_m$-invariant is then defined by
\[
\overline{\mu}_m(I) := \mu_{m}(I) \mod \Delta_m(I).
\]
The following analogue of Theorem 2.1.2 is due to Morishita.
\\
\\
\quad {\bf Theorem 2.2.3} ([Mo3, Theorems 1.2.1, 1.2.5]){\bf .}  
 {\it 
$(1)$ $\zeta_m^{\mu_m(ij)} = \left( \frac{p_j}{p_i} \right)_m$ where $\zeta_m$ is the primitive $m$-th root of unity given in $(2.2.1)$ and $\left( \frac{p_j}{p_i} \right)_m$ is the $m$-th power residue symbol in $\mathbb{Q}_{p_i}$.\\
$(2)$ If $2 \le |I| \le  2^{e_S}$, $\overline{\mu}_m(I)$ is an invariant depending only on $S$. \\
$(3)$ Let $r$ be an integer such that $2 \le r \le  2^{e_S}$.
For multi-indices $I, J$ such that $|I| + |J| = r - 1$, we have, for any $1 \le i\le n$,
\[
\sum_{H \in {\rm PSh}(I, J) } \overline{\mu}_m(Hi) \equiv 0 \mod \mathrm{g.c.d} \{ \Delta(Hi) \mid H \in {\rm PSh}(I, J) \}.
\]
}
\\
\\
\quad {\bf Example 2.2.4.}  Let $S = \{ p_1, p_2, p_3 \} $ be a triple of distinct prime numbers satisfying the condition (1.2.1) and let $m = 2$.
Then $\mu_2(I) = 0 $ if $|I| \le 2$ and hence, for $|I| = 3$, $\Delta_2(I) = 0$ and $\overline{\mu}_2(I) = \mu_2(I) \in \mathbb{Z}/2\mathbb{Z}$.
The following theorem interprets the R\'{e}dei triple symbol as a Milnor invariant.
\\
\\
\quad {\bf Theorem 2.2.4.1} ([Mo2, Theorem 3.2.5]){\bf .}  
{\it 
Under the above assumption on $\{ p_1, p_2, p_3 \} $ we have
\[
[p_1, p_2, p_3 ] = (-1)^{\mu_{2}(123)}.
\]
}
\\
For example, D. Vogel ([V1, Example 3.14]) showed that for $S = \{ 13, 61, 937 \} $ $\mu_{2}(I) = 0$ $(|I| \le 2)$, $\mu_{2}(I) = 1$ ($I$ is a permutation of 123), $\mu_2(ijk) = 0$ (otherwise).
In view of Example 2.1.3, this triple of prime numbers may be called the {\it Borromean primes}.
\begin{center}
\includegraphics[height=4cm]{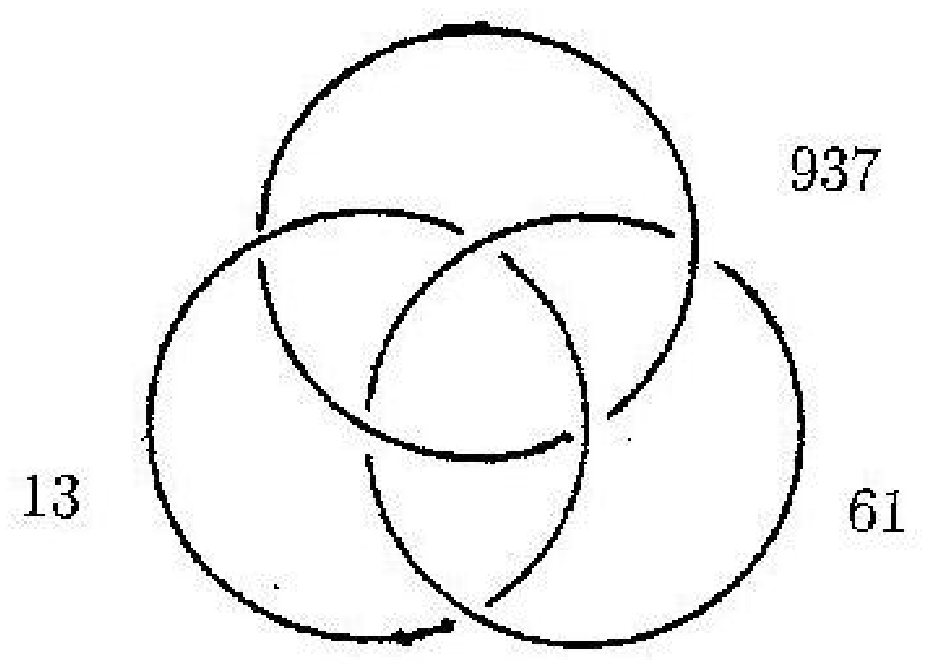}
\end{center}
\quad Finally, we give an analogue of Theorem 2.1.6 for prime numbers.
Let $I = (i_1 \cdots i_n), 2 \le n \le l^{e_S}$ and assume $\Delta_m(I) \neq \mathbb{Z}/m\mathbb{Z}$.
We define the map $\rho_{(m, I)} : \hat{F} \rightarrow N_{n}((\mathbb{Z}/m\mathbb{Z})/\Delta_m(I))$ by
\[
\rho_{(m, I)}(f) := \left(
\begin{array}{ccccc}
1 & \epsilon(\frac{\partial f}{\partial x_{i_1}})_m & \epsilon(\frac{\partial^{2} f}{\partial x_{i_1} \partial x_{i_2}})_m & \cdots & \epsilon(\frac{\partial^{n - 1} f}{\partial x_{i_1} \cdots \partial{x_{i_{n-1}} } })_m\\
& 1 & \epsilon(\frac{\partial f}{\partial x_{i_2} })_m & \cdots &\epsilon(\frac{\partial^{n - 2} f}{\partial x_{i_2} \cdots \partial{x_{i_{n-1}} } })_m\\
&& \ddots & \ddots & \vdots \\
& 0 & & 1 & \epsilon(\frac{\partial f}{\partial x_{i_{n-1}} })_m\\
&&&& 1 \\
\end{array}
\right)
\mod \Delta_m(I),
\] 
where we set $\epsilon(\alpha)_m =  \epsilon_{\mathbb{Z}[[\hat{F}]]}(\alpha)$ mod $m$ for $\alpha \in \mathbb{Z}_{l}[[\hat{F}(l)]]$. It can be shown by the property of the pro-2 Fox derivative that $\rho_{(m,I)}$ is a homomorphism.
\\
\\
\quad {\bf Theorem 2.2.5} ([Mo3, Theorem 1.2.7]){\bf .}  {\it 
$(1)$ The homomorphism $\rho_{(m, I)}$ factors through the Galois group $G_S(2)$.
Further it is surjective if $i_1, \cdots, i_{n - 1}$ are all distinct.
\\
$(2)$ Suppose that $i_1, \ldots, i_{n - 1}$ are all distinct.
Let $K_{(m,I)}$ be the extension over $\mathbb{Q}$ corresponding to $\mathrm{Ker}(\rho_{(m, I)})$.
Then $K_{(m,I)}/\mathbb{Q}$ is a Galois extension unramified outside $p_{i_1}, \dots, p_{i_{n - 1}}$ and $\infty$ with Galois group $\mathrm{Gal}(K_{(m,I)}/\mathbb{Q}) = N_{n}((\mathbb{Z}/m\mathbb{Z})/\Delta_m(I))$.
For a Frobenius automorphism $\sigma_{i_n}$ over $p_{i_n}$, one has 
\[
\rho_{(m, I)}(\sigma_{i_n}) = \left(
\begin{array}{ccccc}
1 & 0  & \cdots & 0 & \overline{\mu}_m(I)  \\ 
& 1 & \cdots &  & 0 \\
&& \ddots &  & \vdots \\
& 0 & & 1 & 0 \\
&&&& 1 \\
\end{array}
\right)
\]
and hence the following holds$:$
\[
\overline{\mu}_m(I) = 0 \Longleftrightarrow p_{i_n} is\ completely\ decomposed\ in\ K_{(m,I)}/\mathbb{Q}.
\]
}
\\
\quad {\bf Example 2.2.6.}  Let $m = 2$ and $K = K_{(2,I)}$. For $S = \{ p_1, p_2 \} $, $p_i \equiv 1 \pmod{4}$ $(i = 1,2)$ and $I=(12)$, we have
$$K = \mathbb{Q}(\sqrt{p_1}), \; {\rm Gal}(K/\mathbb{Q}) = N_2(\mathbb{F}_2) = \mathbb{Z}/2\mathbb{Z}, \; (-1)^{\mu_2(12)} = \left( \frac{p_1}{p_2} \right).$$
For $S = \{ p_1, p_2, p_3 \} $ satisfying the condition (1.2.1) and $I=(123)$, we have 
$$K = k_{\{p_1, p_2\}}, \; {\rm Gal}(K/\mathbb{Q}) = N_3(\mathbb{F}_2) = D_8, \; (-1)^{\mu_2(123)} = [p_1, p_2, p_3].$$
\vspace{.2cm}

 Theorem 2.2.5 suggests a problem to construct concretely a Galois extension $K_n/\mathbb{Q}$ unramified outside $p_1,\dots, p_{n-1}$ and $\infty$ with Galois group $N_n(\mathbb{F}_2)$  and to introduce the multiple residue symbol $[p_1,\dots , p_n]$, as a generalization of the Legendre symbol and the R\'{e}dei triple symbol, which should describe the decomposition law of $p_n$ in the extension $K_n/\mathbb{Q}$ and coincide with  $(-1)^{\mu_{2}(12 \cdots n)}$. In the next section, we solve this problem for the case $n = 4$.
\\

\begin{center}
{\bf \S 3 \quad Construction of an $N_4(\mathbb{F}_2)$-extension and the 4-th multiple residue symbol}
\end{center}

 In this section, under certain conditions on three prime numbers $p_1, p_2, p_3,$ we construct concretely a Galois extension $K$ over $\mathbb{Q}$ where all ramified prime numbers are $p_1, p_2$ and $p_3$ and the Galois group is $N_4(\mathbb{F}_2)$, and introduce the 4-th multiple residue symbol $[p_1, p_2, p_3, p_4]$ which describes the decomposition law of $p_4$ in $K/\mathbb{Q}$. We then show that $[p_1, p_2, p_3, p_4]$ coincides with $(-1)^{\mu_{2}(1234)}$, where $\mu_{2}(1234)$ is the 4-th arithmetic Milnor invariant defined in 2.2. We keep the same notations as in the previous sections.
\\

\quad {\bf3.1. Construction of an $N_4(\mathbb{F}_2)$-extension.}  Let $p_1, p_2$ and $p_3$ be three prime numbers satisfying the conditions 
$$ \left\{ \begin{array}{l}
p_i \equiv 1 \pmod{4} \; (i = 1, 2, 3), \; \displaystyle{\left( \frac{p_i}{p_j} \right) = 1} \; (1 \le i \neq j \le 3), \\
  \mbox{ $[p_i , p_j , p_k]$ } = 1 \; ( \{ i, j, k \} = \{1, 2, 3 \} ).  
\end{array}
\right.  \leqno{(3.1.1)}
$$

We let 
$$ \left\{ \begin{array}{l}
k_i := \mathbb{Q}(\sqrt{p_i})\; (i = 1,2,3), \; k_{ij} := k_ik_j = \mathbb{Q}(\sqrt{p_i},\sqrt{p_j})\; (1\leq i<j\leq 3),\\
 k_{123} := k_1k_2k_3 = \mathbb{Q}(\sqrt{p_1},\sqrt{p_2},\sqrt{p_3}).
\end{array} \right.
$$
For simplicity, we set $k := k_1$ in the following. Let $\frak{p}_2$ be one of prime ideals of ${\cal{O}}_k$ lying over $p_2$.
Then as in Lemma 1.1.2, we can find a triple of integers $(x, y, z)$ with $\alpha = x + y \sqrt{p_1}$ satisfying (1), (2) in Lemma 1.1.2 such that
\[
(\alpha) = \frak{p}_2^{m} \ (\mbox{$m$ being an odd integer}),\quad
k_{ \{p_1,p_2 \} } = \mathbb{Q}(\sqrt{p_1},\sqrt{p_2},\sqrt{\alpha}). 
\]

In the following, we fix such an $\alpha$ once and for all. 
\\

For a prime $\frak{p}$ of $k$, we denote by $\left( \frac{\ ,\ }{ \frak{p} } \right)$ the Hilbert symbol  in the local field $k_{\frak{p}}$,  namely, 
\[
(a, k_{\frak{p}}(\sqrt{b})/k_{\frak{p}})\sqrt{b} = \left( \frac{a, b}{\frak{p}} \right) \sqrt{b} \quad (a, b \in k_{\frak{p}}^{\times}),
\]
where $(\ , k_{\frak{p}}(\sqrt{b})/k_{\frak{p}}) : k_{\frak{p}}^{\times}  \rightarrow \mbox{Gal}(k_{\frak{p}}(\sqrt{b})/k_{\frak{p}})$ is the norm residue symbol of local class field theory. 
\\
\\
\quad {\bf Lemma 3.1.2.}  
 {\it 
For any prime $\frak{p}$ of $k$, we have
\[
\left( \frac{\alpha ,p_3 }{ \frak{p} } \right) = 1.
\]
}
\begin{proof}  We consider the following five cases.\\
(case 1) $\frak{p}$ is prime to $\frak{p}_2$, $p_3$, 2, $\infty$:  Then we have $\alpha ,p_3 \in U_{\frak{p}}$, where $U_{\frak{p}}$ is the unit group of $k_{\frak{p}}$, and hence $\left( \frac{\alpha ,p_3 }{ \frak{p} } \right) = 1$. \\
(case 2) $\frak{p} = \frak{p}_2$:  
Let $\pi$ be a prime element of $k_{\frak{p}_2}$.
Write $\alpha = u_1 \pi^{m_2},$ $u_1 \in U_{\frak{p}_2}$.
Then we have
\begin{align*}
\left(\frac{\alpha, p_3}{\frak{p}_2}\right) 
&= \left(\frac{u_1, p_3}{\frak{p}_2}\right) \left(\frac{\pi^{m_2}, p_3}{\frak{p}_2}\right) \\
&= \left(\frac{\pi, p_3}{\frak{p}_2}\right) \quad (u_1, p_3 \in U_{\frak{p}_2}, m_2 \ \mbox{is odd}) \\
&= \frac{ (\pi, k_{ \frak{p}_2 }(\sqrt{p_3})/k_{ \frak{p}_2 } )\sqrt{p_3} } {\sqrt{p_3}}.
\end{align*}
Since $(\pi, k_{ \frak{p}_2 }(\sqrt{p_3})/k_{ \frak{p}_2 } )$ is the Frobenius automorphism over $p_2$ in $k(\sqrt{p_3})/k$, $(\pi, k_{ \frak{p}_2 }(\sqrt{p_3})/k_{ \frak{p}_2 } )(\sqrt{p_3}) = \sqrt{p_3}$ by $\left( \frac{p_1}{p_2} \right) = \left( \frac{p_3}{p_2} \right) = 1$. \\
(case 3) $\frak{p} \mid p_3$:  Let $\varpi$ be a prime element of $k_\frak{p}$.
Write $p_3 = u_2 \varpi$, $u_2 \in U_{\frak{p}}$.
Then we have
\begin{align*}
\left(\frac{\alpha, p_3}{\frak{p}}\right)
&= \left(\frac{p_3, \alpha}{\frak{p}}\right) \\
&= \left(\frac{u_2, \alpha}{\frak{p}}\right) \left(\frac{\varpi, \alpha}{\frak{p}}\right) \\
&= \left(\frac{\varpi, \alpha}{\frak{p}}\right) \quad (u_2, \alpha \in U_{\frak{p}}) \\
&= \frac{ (\varpi, k_{ \frak{p} }(\sqrt{\alpha})/k_{ \frak{p} } )\sqrt{\alpha} } {\sqrt{\alpha}}.
\end{align*}
Since $\frak{p}$ is decomposed in $k(\sqrt{\alpha})/k$ by $[p_1, p_2, p_3] = 1$ and $(\varpi, k_{ \frak{p} }(\sqrt{\alpha})/k_{ \frak{p} } )$ is the
Frobenius automorphism over $\frak{p}$ in $k(\sqrt{\alpha})/k$, $(\varpi, k_{ \frak{p} }( \sqrt{\alpha} )/k_{ \frak{p} } )( \sqrt{\alpha} ) = \sqrt{\alpha}$. \\
(case 4) $\frak{p} = \infty$:  Since $p_3 > 0$, $\left(\frac{\alpha, p_3}{\infty}\right) = 1$. \\
(case 5) $\frak{p} \mid 2$:  If $\frak{p} = (2)$, the above cases and the product formula for the Hilbert symbol yields $\left(\frac{\alpha, p_3}{\frak{p}}\right) = 1$. 
If $(2) = \frak{p} \cdot \frak{p}'$ $(\frak{p} \neq \frak{p}')$, $k_{ \frak{p} } = k_{ \frak{p}' } = \mathbb{Q}_2$ and so we have
\[
\left(\frac{\alpha, p_3}{\frak{p}}\right) = \left(\frac{\alpha, p_3}{\frak{p}'}\right) = (-1)^{ \frac{p_3 -1}{2} \cdot \frac{\alpha -1}{2} } = 1. \qedhere 
\]
\end{proof} $\;$
\\
\quad {\bf Proposition 3.1.3.}  {\it 
Assume that the class number of $k$ is $1$.
Then there are $X,Y,Z \in {\cal O}_k$ satisfying the following conditions$:$ \\
$(1)$ \quad $X^2 - p_3 Y^2 - \alpha Z^2 = 0 $, \\ 
$(2)$ \quad $\mathrm{g.c.d}(X, Y, Z) = 1$.
}
\begin{proof}  By Lemma 3.1.2, we have $\alpha \in N_{k_{\frak{p}}(\sqrt{p_3})/k_{\frak{p}} }(k_{\frak{p}}(\sqrt{p_3})^{\times})$ for any prime $\frak{p}$ of $k$ and so there are $X_{\frak{p}}, Y_{\frak{p}} \in k_{\frak{p}}$ such that  $X_{\frak{p}}^2 - p_3 Y_{\frak{p}} ^2 = \alpha$.
By the Hasse principal, there are $\tilde{X}, \tilde{Y} \in k$ such that $\tilde{X}^2 - p_3 \tilde{Y}^2 = \alpha$ from which the condition (1) holds by writing $\tilde{X} = \frac X Z, \tilde{Y} = \frac Y Z$ with $X,Y,Z \in {\cal O}_k$.
Since ${\cal O}_k$ is the principal ideal domain by the assumption, we may choose $X,Y,Z \in {\cal O}_k$ so that the condition (2) is satisfied.
\end{proof}

For $k_{13} = \mathbb{Q}(\sqrt{p_1},\sqrt{p_3})$,
let $U$ be the unit group of ${\cal O}_{k_{13}}/(4)$ and $U(2)$ the $2$-Sylow subgroup of $U$. 
Similarly, let $k_{13}' := \mathbb{Q}(\sqrt{p_1},\sqrt{\alpha})$ and define $U':= ({\cal O}_{k_{13}'}/(4))^{\times}$ and $U'(2)$ to be the $2$-Sylow subgroup of $U'$.\\
\\
\quad {\bf Lemma 3.1.4.}  {\it  The group $U(2)$ is given by
\begin{align*}
U(2) &= \langle -1 \rangle \times \langle \sqrt{p_1} \rangle \times \langle \sqrt{p_3} \rangle \times \left\langle \frac{3 + \sqrt{p_1} + \sqrt{p_3} + \sqrt{p_1 p_3}}{2} \right\rangle   \\
	&\simeq \mathbb{Z}/2\mathbb{Z} \times \mathbb{Z}/2\mathbb{Z} \times \mathbb{Z}/2\mathbb{Z} \times \mathbb{Z}/2\mathbb{Z}.
\end{align*}
Similarly, $U'(2)$ is given by
\begin{align*}
U'(2) &= \langle -1 \rangle \times \langle \sqrt{p_1} \rangle \times \langle \sqrt{\alpha} \rangle \times \left\langle \frac{3 + \sqrt{p_1} + \sqrt{\alpha} + \sqrt{p_1 \alpha}}{2} \right\rangle   \\
	&\simeq \mathbb{Z}/2\mathbb{Z} \times \mathbb{Z}/2\mathbb{Z} \times \mathbb{Z}/2\mathbb{Z} \times \mathbb{Z}/2\mathbb{Z}.
\end{align*}}
\begin{proof}  Since 2 is unramified in the extension $k_{13}/\mathbb{Q}$, we have the decomposition $(2) = \frak{c_1} \cdots \frak{c}_r$. Therefore the order of $U$ is given by 
$$
\Pi^{r}_{i=1} N\frak{c_i} (N\frak{c_i} - 1) = N( (2) ) \Pi^{r}_{i=1}(N\frak{c_i} - 1) = 16m 
$$
and  so $U$ has the order $16m$, where $m$ is an odd integer. Let $A := \{ \pm1\mod (4), \pm \sqrt{p_1} \mod (4), \pm \sqrt{p_3} \mod (4), \pm \sqrt{p_1 p_3} \mod (4) \}$. Since $p_i \equiv 1 \mod 4$, each element of $A$ has the order $2$ and so $A \subset U(2)$. We show that the order of $A$ is $8$. Suppose $\sqrt{p_1} \equiv \sqrt{p_3} \mod (4)$ for example.
Then $\sqrt{p_1} - \sqrt{p_3} = 4 \beta$ for some $\beta \in {\cal O}_{k_{13}}$.
Taking the norm $N_{ k_{13} /\mathbb{Q} }$, we obtain
$$
\frac{ -p_1 - p_3}{16} + \frac{ \sqrt{p_1 p_3} }{8} \in {\cal O}_{\mathbb{Q}(\sqrt{p_1 p_3}) } = \{ \frac{a + b \sqrt{p_1 p_3} }{2} | \; a,b \in \mathbb{Z}, a \equiv b \mod 2 \},
$$
which is a contradiction.
Similary, using the structure of ${\cal O}_{k_1}$ and ${\cal O}_{k_3}$, we can check that any two elements in $A$ are distinct. Hence we see that $U(2)  = A \cup A \cdot \{ (3 + \sqrt{p_1} + \sqrt{p_3} + \sqrt{p_1 p_3})/2 \mod (4) \}$. 
Replacing  $p_3$ by $\alpha$, the assertion for $U'(2)$ can be shown similarly.
\end{proof} $\;$
\\
\quad {\bf Lemma 3.1.5.}  {\it 
Assume $p_1 \equiv 5 \pmod{8}$.
Then there is a unit $\epsilon \in {\cal{O}}_{k}^{\times}$ of the form $\epsilon = s + t \sqrt{p_1}, s,t \in \mathbb{Z}, s \equiv 0, t \equiv 1 \pmod{2}$.
Such a unit $\epsilon$ satisfies $\epsilon \equiv \pm \sqrt{p_1} \mod(4)$ in $U(2)$ and $U'(2)$.
}
\begin{proof}  Since $p_1 \equiv 1 \pmod{4}$, the fundamental unit $\epsilon_1 = \frac{s_1 + t_1 \sqrt{p_1} }{2}$ $(s_1 \equiv t_1 \pmod{2})$ of $k$ satisfies $N_{k/\mathbb{Q}}(\epsilon_1) = -1$.
If $s_1 \equiv t_1 \equiv 0 \pmod{2}$, we let $\epsilon := \epsilon_1 = s + t \sqrt{p_1}$, $s := s_1/2, t := t_1/2 \in \mathbb{Z}$, where we have $s \equiv 0, t \equiv 1 \pmod{2}$, since $s^2 -p_1 t^2 = -1$.
Since $\epsilon = s + t \sqrt{p_1} = s + s \sqrt{p_1} + (t - s)\sqrt{p_1}$ and $s + s \sqrt{p_1} \in 4{\cal{O}}_{k_{13}}$, $\epsilon 
\equiv \pm \sqrt{p_1} \mod(4)$.
Suppose $s_1 \equiv t_1 \equiv 1 \pmod{2}$.
Since $p_1 \equiv 5 \pmod{8}$, we have $s_1^2 + 3 p_1 t_1^2 \equiv 3 s_1^2 + p_1 t_1^2 \equiv 0 \pmod{8}$ and so
\[
\epsilon_1^3 = \frac{s_1(s_1^2 + 3 p_1 t_1^2) + t_1(3 s_1^2 + p_1t_1^2)\sqrt{p_1} }{8} = s + t \sqrt{p_1},
\]
where $s = s_1(s_1^2 + 3 p_1 t_1^2)/8, t = t_1(3 s_1^2 + t_1^2)/8 \in \mathbb{Z}$. Since $N_{k/\mathbb{Q}}(\epsilon_1^3) = -1$, $\epsilon = \epsilon_1^3$ satisfies the desired conditions. 
\end{proof}
The  following theorem may be regarded as an analogue of Lemma 1.1.2.
\\
\\
\quad {\bf Theorem 3.1.6.} {\it 
Assume that the class number of $k$ is $1$ and $p_1 \equiv 5 \pmod{8}$.
Then there are $X,Y,Z \in {\cal O}_k$ satisfying the following conditions$:$\\ 
$(1)$ \quad $X^2 - p_3 Y^2 - \alpha Z^2 = 0 $, \\ 
$(2)$ \quad $\mathrm{g.c.d}(X, Y, Z) = 1$, $(Z,2) = 1$ $(\mathrm{resp. \ g.c.d}(X, Y, Z) = 1$, $(Y,2) = 1)$,
$(3)$ \quad There is $\lambda \in {\cal O}_{k_{13}}$ $(\mathrm{resp.} \ \lambda \in {\cal O}_{k_{13}'})$ such that $\lambda^2 \equiv X + Y \sqrt{p_3} \mod (4)$  $(\mathrm{resp.} \ \lambda^2 \equiv X + Z \sqrt{\alpha} \mod (4))$.
}
\begin{proof} By Proposition 3.1.3, there are $X,Y,Z \in {\cal O}_k$ satisfying (1) and (2). 

Case $(Z,2) = 1$: \quad 
Let $\theta := X + Y\sqrt{p_3}$ and $\overline{\theta} := \theta \mod (4)$.
Then we easily see $\theta \in {\cal{O}}_{k_{13}}$ and $\overline{\theta} \in U$ since $(Z,2) = 1$. Let $n$ be the order of $\overline{\theta}$ in $U$. 
\\
(i) Suppose $n \not\equiv 0 \pmod{2}$.
Then it is easy to see that there is $\lambda \in {\cal O}_{k_{13}}$ such that $\lambda^2 \equiv \theta \mod (4)$.  
\\
(ii) Suppose $n \equiv 0 \pmod{2}$.
By Lemma 3.1.4, $\frac{n}{2} \not\equiv 0 \pmod{2}$ and $ \overline{\theta}^{\frac{n}{2}} \in U(2)$.
Write $\theta^{\frac{n}{2}} = b_1 + b_2 \sqrt{p_1} + b_3 \sqrt{p_3} + b_4 \sqrt{p_1 p_3} $, $b_i \in \mathbb{Q}$.
Since $N_{k_{13}/k}(\theta) = X^2 - p_3 Y^2 = \alpha Z^2 $, $N_{k_{13}/k}(\theta^{\frac{n}{2}}) = (\alpha Z^2 )^{\frac{n}{2}}$.
Since $\alpha = x + y \sqrt{p_1} = x - y + 2y \cdot \frac{1+\sqrt{p_1}}{2} \equiv 1 \mod (4)$, $(\alpha Z^2 )^{\frac{n}{2}} \equiv (Z^{\frac{n}{2}} )^2 \equiv 1\mod (4)$.
Therefore we have 
\[
(b_1 + b_2 \sqrt{p_1} + b_3 \sqrt{p_3} + b_4 \sqrt{p_1 p_3} ) \cdot (b_1 + b_2 \sqrt{p_1} - b_3 \sqrt{p_3} - b_4 \sqrt{p_1 p_3} ) \equiv 1\mod (4). \leqno{(3.1.6.1)}
\]
We claim that $\theta^{\frac{n}{2}} \equiv -1 \  \mbox{or} \  \pm \sqrt{p_1} \mod (4)$. 
Suppose this is not the case.
Then, by Lemma 3.1.4, $\theta^{\frac{n}{2}} \equiv \pm \sqrt{p_3}, \  \pm \sqrt{p_1 p_3} \  \mbox{or} \  a \cdot (3 + \sqrt{p_1} + \sqrt{p_3} + \sqrt{p_1 p_3})/2 \ (a \in A) \mod (4)$ and so the coefficients of $\sqrt{p_3}$ or $\sqrt{p_1 p_3} $ are not $0$.
Since any element of $U(2)$ has order $2$, we have 
\[
(b_1 + b_2 \sqrt{p_1} + b_3 \sqrt{p_3} + b_4 \sqrt{p_1 p_3} ) \cdot (b_1 + b_2 \sqrt{p_1} - b_3 \sqrt{p_3} - b_4 \sqrt{p_1 p_3} ) \not\equiv 1 \mod (4),
\]
which contradicts to (3.1.6.1).
Therefore, by Lemma 3.1.5, there is $\epsilon \in {\cal{O}}_{k}^{\times}$ such that $\epsilon \theta^{\frac{n}{2}} \equiv 1 \mod (4)$ and $\epsilon^2 \equiv 1 \mod (4)$.
Replacing $(X,Y,Z)$ by $(\epsilon X,\epsilon Y,\epsilon Z)$, (1), (2) holds obviously, and (3) is also satisfied because $\epsilon \theta \mod (4)$ has the order $\frac{n}{2} \not\equiv 0 \pmod{2}$. 

Case $(Y,2) = 1$:  Let $\theta' := X + Z \sqrt{\alpha}$. Replacing $\theta$ by $\theta'$ and $p_3$ by $\alpha$, the above proof works well by using Lemma 3.1.4.\end{proof}

Let $\boldsymbol{a} = (X,Y, Z)$ be a triple of integers in ${\cal O}_{k}$ satisfying (1), (2), (3) in Theorem 3.1.6 and fix it once and for all.
We let
$$ \left\{ \begin{array}{ll}
\theta := X + Y \sqrt{p_3} &  \; \mbox{if}\;  (Z,2) = 1,\\
\theta' := X + Z \sqrt{\alpha} &  \; \mbox{if}\;  (Y,2) = 1,
\end{array} \right.
$$
and set
$$ \left\{ 
\begin{array}{l} \theta_1 := \theta,\\
\theta_2 := X - Y \sqrt{p_3}, \\
\theta_3 := \overline{X} + \overline{Y} \sqrt{p_3},\\
\theta_4 := \overline{X} - \overline{Y} \sqrt{p_3},
\end{array}\right.
\;\; \;\;
 \left\{
\begin{array}{l} 
\theta'_1 = \theta',\\
\theta'_2 = X - Z \sqrt{\alpha}, \\
\theta'_3 = \overline{X} + \overline{Z} \sqrt{\overline{\alpha} },\\ 
\theta'_4 = \overline{X} - \overline{Z} \sqrt{\overline{\alpha}},\\
\end{array}\right.
$$
where $\overline{X}, \overline{Y}$ and $\overline{\alpha}$ are conjugates of $X$, $Y$ and $\alpha$ over $\mathbb{Q}$ respectively.\\
\\
\quad {\bf Definition 3.1.7.} We then define the number field $K$ by
$$
K = K_{ \boldsymbol{a} } = \left\{
\begin{array}{ll}
\mathbb{Q}(\sqrt{p_1}, \sqrt{p_2}, \sqrt{p_3}, \sqrt{\theta_1 \theta_2}, \sqrt{\theta_1 \theta_3}, \sqrt{\theta_1})
 &\quad   \mbox{if}\  (Z,2)  = 1, \\
\mathbb{Q}(\sqrt{p_1}, \sqrt{p_2}, \sqrt{p_3}, \sqrt{\theta_1'\theta_2'}, \sqrt{\theta_1' \theta_3'}, \sqrt{\theta_1'})
&\quad \mbox{if}\  (Y,2)  = 1.
\end{array}
\right.  
$$
\\
For the latter use, we set, for the case of $(Y,2) = 1$, 
$$ \left\{ \begin{array}{l} 
\eta_1 := (\sqrt{\theta'_1} + \sqrt{\theta'_2})^2 = 2X + 2Y \sqrt{p_3}, \\
\eta_2 := (\sqrt{\theta'_1} - \sqrt{\theta'_2})^2 = 2X - 2Y \sqrt{p_3}, \\ 
\eta_3 := (\sqrt{\theta'_3} + \sqrt{\theta'_4})^2 = 2\overline{X} + 2\overline{Y} \sqrt{p_3},\\ 
\eta_4 := (\sqrt{\theta'_3} - \sqrt{\theta'_4})^2 = 2\overline{X} - 2\overline{Y} \sqrt{p_3}.
\end{array}\right.
$$
\\
\quad {\bf Theorem 3.1.8}  {\it 
$(1)$ We have}
$$
K = \left\{
\begin{array}{ll}
\mathbb{Q}(\sqrt{\theta_1}, \sqrt{\theta_2}, \sqrt{\theta_3}, \sqrt{\theta_4})
& \quad   \mathrm{if}\  (Z,2)  = 1, \\
 \mathbb{Q}(\sqrt{\theta_1'}, \sqrt{\theta_2'}, \sqrt{\theta_3'}, \sqrt{\theta_4'}) 
= \mathbb{Q}(\sqrt{\eta_1}, \sqrt{\eta_2}, \sqrt{\eta_3}, \sqrt{\eta_4})&\quad \mathrm{if}\  (Y,2)  = 1.
\end{array}
\right. 
$$
$(2)$ {\it The extension $K/\mathbb{Q}$ is a Galois extension whose Galois group is isomorphic to $N_{4}(\mathbb{F}_2)$.
}
\begin{proof}  (1) Case $(Z,2) = 1$:  It is easy to see $\sqrt{\theta_2}, \sqrt{\theta_3} \in K$. Noting that 
$$
\left.
\begin{array}{ll}
\theta_1 \theta_2 \theta_3 \theta_4 
&= N_ {k_{13}/\mathbb{Q} }( \theta_1 ) \\
&= N_{ k/\mathbb{Q} }( N_{ k_{13}/k }( \theta_1 ) ) \\
&= N_{ k/\mathbb{Q} }(\alpha Z^2) \\
&= p_2 h^2 \quad (h \in \mathbb{Z}),
\end{array}
\right. \leqno{(3.1.8.1)}
$$
 we have $\sqrt{\theta_4} \in K$ and hence $\mathbb{Q}(\sqrt{\theta_1}, \sqrt{\theta_2}, \sqrt{\theta_3}, \sqrt{\theta_4}) \subset K$. Next we show the converse inclusion. Write $\theta_1 = a_1 + a_2 \sqrt{p_1} + a_3 \sqrt{p_3} + a_4 \sqrt{p_1 p_3}$ ($a_i \in \mathbb{Q}$). By considering the prime factorization of the ideal $(\alpha Z^2 )$ in $k_1$, we find $\alpha Z^2 \notin \mathbb{Z}$. Then, by the equality $\theta_1 \theta_2 = \alpha Z^2$, we find that the number of $i$ $(1 \le i \le 4)$ with $a_i = 0$ is at most one. Since $\theta_1 + \theta_2 = 2 (a_1 + a_2 \sqrt{p_1}), \theta_1 +  \theta_3 = 2(a_1 + a_3 \sqrt{p_3})$ and $\theta_1 + \theta_4 = 2(a_1 + a_4 \sqrt{p_1 p_3})$,  $\sqrt{p_1}, \sqrt{p_3} \in \mathbb{Q}(\sqrt{\theta_1}, \sqrt{\theta_2}, \sqrt{\theta_3}, \sqrt{\theta_4})$.
By (3.1.8.1), we get $ K \subset \mathbb{Q}(\sqrt{\theta_1},\sqrt{\theta_2},\sqrt{\theta_3}, $ $\sqrt{\theta_4})$. 

Case $(Y,2) = 1$: First, let us show $ \mathbb{Q}(\sqrt{\theta_1'}, \sqrt{\theta_2'}, \sqrt{\theta_3'}, \sqrt{\theta_4'}) 
= \mathbb{Q}(\sqrt{\eta_1}, \sqrt{\eta_2},$ $\sqrt{\eta_3}, \sqrt{\eta_4})$.
By the definition of $\eta_i$'s, obviously the inclusion $\supset$ holds.
Since $\sqrt{\eta_1} + \sqrt{\eta_2} = 2\sqrt{\theta'_1},$  $\sqrt{\eta_1} - \sqrt{\eta_2} = 2\sqrt{\theta'_2}$, $\sqrt{\eta_3} + \sqrt{\eta_4} = 2\sqrt{\theta'_3}$,
$\sqrt{\eta_3} - \sqrt{\eta_4} = 2\sqrt{\theta'_4}$, we obtain the converse inclusion $\subset$.

Next, we show $K = \mathbb{Q}(\sqrt{\theta_1'}, \sqrt{\theta_2'}, \sqrt{\theta_3'}, $ $\sqrt{\theta_4'})$. It is easy to see  $\sqrt{\theta_2'}, \sqrt{\theta_3'} \in K$. Since $\theta_1'\theta_2' = X^2 - \alpha Z^2 = p_3 Y^2$, we have $\theta_3'\theta_4' = \overline{X}^2 - \overline{\alpha} \overline{Z} ^2 = p_3 \overline{Y}^2$. So,  $\theta_1' \theta_2' \theta_3' \theta_4' = p_3^2(Y \overline{Y})^2 \in \mathbb{Q}$ and $\sqrt{\theta_4'} \in K $.
For the converse inclusion, it suffices to show $K \subset \mathbb{Q}(\sqrt{\eta_1}, \sqrt{\eta_2}, \sqrt{\eta_3}, \sqrt{\eta_4})$.
By considering the prime factorization of the ideal $(\alpha (2Z)^2 )$ in $k_1$, we find $\alpha (2Z)^2 \notin \mathbb{Z}$.
By $N_{k_{13}/\mathbb{Q}}(\eta_1) = 4p_2 h^2$ and the argument similar to the case of $(Z,2)=1$, we have $\sqrt{p_i} \in \mathbb{Q}(\sqrt{\eta_1}, \sqrt{\eta_2}, \sqrt{\eta_3}, \sqrt{\eta_4})$ ($i = 1,2,3$). \\
\quad (2) Case $(Z,2) = 1$: First, $K/\mathbb{Q}$ is a Galois extension, because $K$ is the splitting field of $\Pi^{4}_{i=1}(T^2 - \theta_i) = \Pi_{\sigma \in \mathrm{Gal}(k_{13}/\mathbb{Q})}(T^2 - \sigma(\theta_1)) \in \mathbb{Z}[T]$.
Next, let $k_{123} = \mathbb{Q}(\sqrt{p_1}, \sqrt{p_2}, \sqrt{p_3})$, $K_1 := k_{123}(\sqrt{\theta_1 \theta_2})$ and $K_2 := K_1( \sqrt{\theta_1 \theta_3})$. Since $\theta_3\theta_4 = \overline{\theta_1}\overline{\theta_2}$ and $\sqrt{\theta_3\theta_4} = h\sqrt{p_2}/\sqrt{\theta_1\theta_2} \in K_1$, $K_1/k_{123}$ is a Galois extension.
Let us show $[K_1 : k_{123}] = 2$.
Define $\sigma \in \mbox{Gal}(k_{123} /\mathbb{Q})$ by 
$$\sigma : (\sqrt{p_1}, \sqrt{p_2}, \sqrt{p_3} ) \mapsto  (-\sqrt{p_1}, -\sqrt{p_2}, \sqrt{p_3}).$$
Let $\tilde{\sigma} \in \mbox{Gal}(K_1/\mathbb{Q})$ be an extension of $\sigma$.
Then we have 
$$(\tilde{\sigma}(\sqrt{\theta_1 \theta_2}))^2 = \tilde{\sigma}(\theta_1 \theta_2) = \theta_3\theta_4$$
and so $\tilde{\sigma}(\sqrt{\theta_1 \theta_2}) = \pm \sqrt{\theta_3\theta_4}$.
Therefore we have
$$\tilde{\sigma}^2(\sqrt{\theta_1 \theta_2}) = \tilde{\sigma}(\pm \sqrt{\theta_3\theta_4}) = \tilde{\sigma}(\pm h \sqrt{p_2}/ \sqrt{\theta_1\theta_2}) = -\sqrt{\theta_1\theta_2}.$$
Since $\tilde{\sigma}^2 \mid_{k_{123}} = \mbox{id}$, $\sqrt{\theta_1\theta_2} \not\in k_{123}$ and hence $[K_1 : k_{123}] = 2$.  Similarly we can show that $K_2/K_1$ is a Galois extension and $[K_2 : K_1] = [K : K_2] = 2$.
Hence we have $[K : \mathbb{Q}] = [K : K_2] [K_2 : K_1][K_1 : k_{123}][k_{123} : \mathbb{Q}] = 64$.

Case $(Y,2) = 1$: $K/\mathbb{Q}$ is a Galois extension, because $K$ is the splitting field of $\Pi^{4}_{i=1}(T^2 - \eta_i) = \Pi_{\sigma \in \mathrm{Gal}(k_{13}/\mathbb{Q})}(T^2 - \sigma(\eta_1)) \in \mathbb{Z}[T]$. Let $K'_1 := k_{123}(\sqrt{\eta_1 \eta_2})$ and $E'_2 := k_{123}(\sqrt{\eta_1 \eta_3})$. By the argumet similar to the case $(Z,2) = 1$, we have $[K : \mathbb{Q}] = [K : K'_2] [K'_2 : K'_1][K'_1 : k_{123}][k_{123} : \mathbb{Q}] = 64$.

Finally, by the computer calculation using GAP, we have the following presentation of the group $N_{4}(\mathbb{F}_2)$:
\[
N_4(\mathbb{F}_2) =  \left\langle g_1, g_2, g_3 \left|
\begin{array}{lr} g_1^2 = g_2^2 = g_3^2 = (g_1 g_3)^2 = 1 & \\
 (g_1 g_2)^4 = (g_2 g_3)^4 = (g_1 g_2 g_3)^4 = 1 & \\
 ((g_1 g_2 g_3 g_2)^2 g_3)^2= 1 & \\
\end{array}
\right.
\right\rangle ,
\]
where $g_1, g_2$ and $g_3$ are words representing the following matrices respectively:
\[
g_1 = \left (
\begin{array}{cccc}
1 & 1 & 0 & 0 \\ 
0 & 1 & 0 & 0 \\
0 & 0 & 1 & 0 \\
0 & 0 & 0 & 1 \\
\end{array}
\right), \;\;
g_2 = \left (
\begin{array}{cccc}
1 & 0 & 0 & 0 \\ 
0 & 1 & 1 & 0 \\
0 & 0 & 1 & 0 \\
0 & 0 & 0 & 1 \\
\end{array}
\right),\;\;
g_3 = \left (
\begin{array}{cccc}
1 & 0 & 0 & 0 \\ 
0 & 1 & 0 & 0 \\
0 & 0 & 1 & 1 \\
0 & 0 & 0 & 1 \\
\end{array}
\right).
\]

Case $(Z,2) = 1$: We define $\tau_1, \tau_2 , \tau_3 \in \mbox{Gal}(K/\mathbb{Q})$ by 
\[
\left.
\begin{array}{rl}
\tau_1 : & (\sqrt{p_1},\sqrt{p_2},\sqrt{p_3},\sqrt{\theta_1\theta_2},\sqrt{\theta_1\theta_3},\sqrt{\theta_1}, \sqrt{\theta_2}, \sqrt{\theta_3}, \sqrt{\theta_4}) \\
 &\mapsto (-\sqrt{p_1},\sqrt{p_2},\sqrt{p_3},\sqrt{\theta_3\theta_4},\sqrt{\theta_1\theta_3},\sqrt{\theta_3},\sqrt{\theta_4}, \sqrt{\theta_1}, \sqrt{\theta_2}) \\
\tau_2 : & (\sqrt{p_1},\sqrt{p_2},\sqrt{p_3},\sqrt{\theta_1\theta_2},\sqrt{\theta_1\theta_3},\sqrt{\theta_1}, \sqrt{\theta_2}, \sqrt{\theta_3}, \sqrt{\theta_4}) \\
 &\mapsto (\sqrt{p_1},-\sqrt{p_2},\sqrt{p_3},-\sqrt{\theta_1\theta_2}, -\sqrt{\theta_1\theta_3},-\sqrt{\theta_1},\sqrt{\theta_2}, \sqrt{\theta_3}, \sqrt{\theta_4}) \\
\tau_3 : & (\sqrt{p_1},\sqrt{p_2},\sqrt{p_3},\sqrt{\theta_1\theta_2},\sqrt{\theta_1\theta_3},\sqrt{\theta_1}, \sqrt{\theta_2}, \sqrt{\theta_3}, \sqrt{\theta_4}) \\
 &\mapsto (\sqrt{p_1},\sqrt{p_2},-\sqrt{p_3},\sqrt{\theta_1 \theta_2},\sqrt{\theta_2\theta_4},\sqrt{\theta_2},\sqrt{\theta_1}, \sqrt{\theta_4}, \sqrt{\theta_3}).
\end{array}
\right.
\]
Then we can easily check 
$\tau_1^2 = \tau_2^2 = \tau_3^2 = (\tau_1 \tau_3 )^2 = \mbox{id}$,   $ (\tau_1 \tau_2)^4 = (\tau_2 \tau_3)^4 = (\tau_1\tau_2 \tau_3)^4 = \mbox{id}$, $ ((\tau_1 \tau_2 \tau_3 \tau_2)^2 \tau_3)^2 = \mbox{id}$. Thus the correspondence $\tau_i  \mapsto g_i$ $(i = 1,2,3)$ gives an isomorphism $\mbox{Gal}(K/\mathbb{Q}) \simeq N_{4}(\mathbb{F}_2)$.

Case $(Y,2) = 1$: We note  
$K = \mathbb{Q}(\sqrt{p_1}, \sqrt{p_2}, \sqrt{p_3}, \sqrt{\eta_1 \eta_2}, \sqrt{\eta_1 \eta_3}, \sqrt{\eta_1}),
$ because $\sqrt{\eta_3\eta_4} = 4h\sqrt{p_2}/\sqrt{\eta_1\eta_2} \in K_1$.  Then the assertion can be shown in a way similar to the case $(Z,2)=1$, by replacing $\theta_i$ with $\eta_i$. 
\end{proof} 

Next, let us study the ramification in our extension $K/\mathbb{Q}$.
First, we recall the following well-known fact on the ramification in a Kummer extension.
\\
\\
\quad {\bf Lemma 3.1.9} ([B, Lemma 6]){\bf.}
{\it \quad
Let $l$ be a prime number and $E$ a number field containing a primitive $l$-th root of unity. Let $E(\sqrt[l]{a})$ $(a \in {\cal O}_E)$ be a Kummer extension over $E$ of degree $l$.
Suppose $(a) = \frak{q}^m\frak{a}$ where  $\frak{q}$ is a prime ideal in $E$ which does not divide $l$, $(\frak{q}, \frak{a}) = 1$ and $l \mid m$.
Then $\frak{q}$ is unramified in $E(\sqrt[l]{a})/E$.
}
\\
\\
\quad {\bf Theorem 3.1.10.}
{\it \quad
All prime numbers ramified in the extension $K/\mathbb{Q}$ are $p_1, p_2$ and $p_3$ with ramification index $2$.
}
\\
\begin{proof}  Case $(Z,2) = 1$:  Let us study the ramification in the extension $k_{13}(\sqrt{\theta_1})/k_{13}$. Since  $(T- \frac{\lambda + \sqrt{\theta_1} }{2} ) (T- \frac{\lambda - \sqrt{\theta_1} }{2}) = (T - \frac{\lambda}{2})^2 - (\frac{\sqrt{\theta_1} }{2})^2 = T^2 - \lambda T + \frac{\lambda^2}{4} -\frac{\theta_1}{4}$ with $\lambda, \frac{\lambda^2 -\theta_1 }{4} \in {\cal{O}}_{k_{13}},$ we find $\frac{\lambda + \sqrt{\theta_1} }{2} \in {\cal{O}}_{k_{13}(\sqrt{\theta_1})}$. Since the relative discriminant of $\frac{\lambda+ \sqrt{\theta_1} }{2}$ in $k_{13}(\sqrt{\theta_1})/k_{13}$ is given by 
\[
\begin{aligned} \left|
\begin{array}{cc}
1 & \frac{\lambda + \sqrt{\theta_1} }{2} \\
1 & \frac{\lambda - \sqrt{\theta_1} }{2} \\
\end{array}
\right|^{2} \end{aligned} 
= \left(\frac{\lambda - \sqrt{\theta_1} }{2} - \frac{\lambda + \sqrt{\theta_1} }{2} \right)^{2} = \theta_1,
\]
 we find that any prime factor of $2$ is unramified in $k_{13}(\sqrt{\theta_1})/k_{13}$.

Next, let us look closely at the prime factorization of the ideal $(\theta_1)$ in $k_{13}$. We let
\[
(\theta_1) = \frak{Q}_1^{e_1}\frak{Q}_2^{e_2} \cdots \frak{Q}_r^{e_r}
\]
be the prime factorization of $(\theta_1)$ and let $\frak{q}_i = \frak{Q}_i \cap k$. Since $N_{k_{13}/k}(\theta) = X^2 - p_3 Y^2 = \alpha Z^2 $, we have
\[
N_{ k_{13}/k }( (\theta_1) ) = (\alpha Z^2) = \frak{p}_2^{m} \frak{a}^2,
\leqno{(3.1.10.1)}
\]
where $\frak{a} := (Z)$ is an ideal in $k$. Now the prime factorization of $\frak{q}_i$ in $k_{13}/k$ has the following three cases:

(i)  $\frak{q}_i = \frak{Q}_i^2 \quad N_{ k_{13}/k }( \frak{Q}_i  ) = \frak{q}_i, $ 

(ii) $\frak{q}_i = \frak{Q}_i \quad N_{ k_{13}/k }( \frak{Q}_i  ) = \frak{q}_i^{2} $, 

(iii) $\frak{q}_i = \frak{Q}_i \frak{Q}'_i  \quad N_{ k_{13}/k }( \frak{Q}_i  ) = \frak{q}_i,\ N_{ k_{13}/k }( \frak{Q}'_i  ) = \frak{q}_i $ 
\\
Case (i): If $e_i$ is odd, it contradicts to (3.1.10.1).
 Hence $e_i$ is even. \\
Case (ii): Since $\theta_1 \in \frak{Q}_i $ and $\frak{q}_i = \frak{Q}_i$,  $\theta_2 = a_1 + a_2 \sqrt{p_1} - a_3 \sqrt{p_3} - a_4 \sqrt{p_1 p_3} = X - Y\sqrt{p_3}\in \frak{Q}_i $. Since $\frak{p}_2$ is decomposed in $k_{13}/k$,  we see, by (3.1.10.1), $Z \in \frak{Q}_i $. Further, Since $\frak{Q}_i $ is not a prime factor of 2 by $(Z,2) =1$ and $2X = \theta_1 + \theta_2 \in \frak{Q}_i $, $2Y\sqrt{p_3} = \theta_1 - \theta_2 \in \frak{Q}_i $ and $X, Y, Z\in k$, we have $X, Y, Z\in \frak{q}_i$, which contradicts to g.c.d$(X,Y,Z) = 1$.\\
Case (iii): Suppose $\frak{P}$ and $\frak{P}'$ are prime factors of $\frak{p}_2$. Since the exponent $m$ in (3.1.10.1) is odd, one of $\frak{P}$ and $\frak{P}'$ appears odd times in the prime factorization of $(\theta_1)$. Let $\frak{P}$ be that one. When $\frak{Q}_i \neq \frak{P}$, assume $e_i$ is odd. By (3.1.10.1), $\frak{Q}'_i$ also appears odd times in the prime factorization of $(\theta_1)$. Therefore we have $\theta_1 \in \frak{Q}_i \frak{Q}'_i  = \frak{q}_i$ and  $\theta_2 \in \frak{q}_i$, and so  $2X = \theta_1 + \theta_2 \in \frak{Q}_i $, $2Y\sqrt{p_3} = \theta_1 - \theta_2 \in \frak{Q}_i $. This deduces  $X, Y, Z \in \frak{q}_i$, which contradicts to g.c.d $(X, Y, Z) = 1$. Thus $e_i$ must be even. 

Getting all together, we find that $(\theta_1)$ has the form $\frak{P}^{m_1}\frak{A}^2$ ($m_1$: odd). Then, by Lemma 3.1.9, ramified finite primes in $k_{13}(\sqrt{\theta_1})/k_{13}$ must be lying over $p_2$. Similarly, we see that ramified finite primes in $k_{13}(\sqrt{\theta_i})/k_{13}$ $(i =2, 3, 4)$ are all lying over $p_2$. This shows that any ramified finite prime in the extension $K = k_{13}(\sqrt{\theta_1}, \sqrt{\theta_2}, \sqrt{\theta_3}, \sqrt{\theta_4})/k_{13}$ is lying over $p_2$. Since $k_{13}/\mathbb{Q}$ is unramified outside $p_1, p_3$, we conclude that all ramified prime numbers in $K/\mathbb{Q}$ are $p_1, p_2$ and $p_3$.

Finally, we show that the ramification indices of $p_i$'s in $K/\mathbb{Q}$ are all 2. 
We easily see that this is true for $p_1$ and $p_3$, because the ramification indices of $p_1$ and $p_2$ in $k_{13}/\mathbb{Q}$ are 2 and any prime factor of $p_1$ or $p_3$ is unramified in $K/k_{13}$. So it suffices to show our assertion for $p_2$. Let $\frak{p}_{2i}$ be a prime factor in $k_{13}$ of $p_2$ which is ramified in $k_{13}(\sqrt{\theta_1})/k_{13}$. Since we have $\frak{p}_{2i} = \frak{Q}_i^2 \quad \mbox{in}\; \quad k_{13}(\sqrt{\theta_1})$, by considering the prime factorization of the ideal $(\theta_i)$ in $k_{13}(\sqrt{\theta_1})$, we see by Lemma 3.1.9 that $\frak{Q}_i$ is unramified in $k_{13}(\sqrt{\theta_1}, \sqrt{\theta_i})$. Therefore any prime factor of $p_2$ ramified in $k_{13}(\sqrt{\theta_1})/k_{13}$ is unramified in 
$k_{13}(\sqrt{\theta_1}, \sqrt{\theta_2}, \sqrt{\theta_3}, \sqrt{\theta_4})/k_{13}(\sqrt{\theta_1})$. Thus the ramification index of $p_2$ is 2.

Case $(Y,2)=1$: As in the case of $(Z,2)=1$, we consider the prime factorization of  $(\theta_1')$ in $k_{13}'$. Then, by a similar argument, we find that $(\theta_1')$ has the ideal decomposition of the form $\frak{Q}'\frak{B}^2$ where any prime factor of $\frak{Q}'$ is lying over $p_3$. This shows by Lemma 3.1.9 that  any ramified finite prime in $k_{13}'(\sqrt{\theta_1'})/k_{13}'$ is lying over $p_3$. Similarly, we see that finite ramified primes in $k_{13}'(\sqrt{\theta_2'})/k_{13}'$, $\overline{k_{13}'}(\sqrt{\theta_3'}
)/\overline{k_{13}'}$ and $\overline{k_{13}'}(\sqrt{\theta_4'})/\overline{k
_{13}'}$ are all lying over $p_3$. Hence all ramified prime numbers in 
$K/\mathbb{Q}$ are $p_1, p_2$ and $p_3$. The assertion on the ramification indices of $p_i$'s can also be shown by an argument similar to the case of $(Z,2)=1$.
\end{proof} $\;$
\\
\quad {\bf Theorem 3.1.11.} 
{\it We have}
\[
K = \left\{
\begin{array}{ll}
k_{ \{ p_1, p_2 \} }k_{ \{ p_2, p_3 \} }(\sqrt{\theta_1})
& \quad   \mbox{if}\  (Z,2)  = 1, \\
k_{ \{ p_1, p_2 \} }k_{ \{ p_3, p_2 \} }(\sqrt{\theta_1'}) 
& \quad \mbox{if}\  (Y,2)  = 1.
\end{array}
\right.
\]
\begin{proof}  Case $(Z,2)  = 1$:  First we have 

\begin{align*}
\mathbb{Q}(\sqrt{p_1}, \sqrt{p_2}, \sqrt{\theta_1 \theta_2}) 
&= \mathbb{Q}(\sqrt{p_1}, \sqrt{p_2}, \sqrt{\alpha Z^2}) \\
&= \mathbb{Q}(\sqrt{p_1}, \sqrt{p_2}, \sqrt{\alpha}) \\
&= k_{ \{ p_1, p_2 \} }. 
\end{align*}
Next, it is easy to see that $\mathbb{Q}(\sqrt{p_2}, \sqrt{p_3}, \sqrt{\theta_1 \theta_3})$ is a dihedral extension over $\mathbb{Q}$ of degree 8. Since all prime numbers ramified in $\mathbb{Q}(\sqrt{p_2}, \sqrt{p_3}, \sqrt{\theta_1 \theta_3})/\mathbb{Q}$ are $p_2$ and $p_3$ with ramification index 2 by Theorem 3.1.10, we have 
$$\mathbb{Q}(\sqrt{p_2}, \sqrt{p_3}, \sqrt{\theta_1 \theta_3}) =  k_{ \{ p_3, p_2 \}}$$
by Theorem 1.1.7. Hence we have
\[
K = k_{ \{ p_1, p_2 \} }k_{ \{ p_3, p_2 \} }(\sqrt{\theta_1}).
\]

 Case $(Y,2)  = 1$:  Noting that $\eta_1  = 2X + 2Y \sqrt{p_3}$, 
$\eta_2 =  2X - 2Y \sqrt{p_3}$ and $\eta_3 = 2\overline{X} + 2\overline{Y} \sqrt{p_3}$, we have
\begin{align*}
\mathbb{Q}(\sqrt{p_1}, \sqrt{p_2}, \sqrt{\eta_1 \eta_2}) 
&= \mathbb{Q}(\sqrt{p_1}, \sqrt{p_2}, \sqrt{4 \alpha Z^2}) \\
&= \mathbb{Q}(\sqrt{p_1}, \sqrt{p_2}, \sqrt{\alpha}) \\
&= k_{ \{ p_1, p_2 \}. } 
\end{align*}
By the same argument as in the case of $(Z,2)=1$ replacing $\theta_i$ with $\eta_i$, we have  $\mathbb{Q}(\sqrt{p_2}, \sqrt{p_3}, \sqrt{\eta_1 \eta_3}) 
= k_{ \{ p_3, p_2 \}}$. Hence  we have, by Theorem 3.1.8, 
\begin{align*}
K
&= \mathbb{Q}( \sqrt{\eta_1}, \sqrt{\eta_2}, \sqrt{\eta_3}, \sqrt{\eta_4}) \\
&= \mathbb{Q}(\sqrt{p_1}, \sqrt{p_2}, \sqrt{p_3}, \sqrt{\eta_1 \eta_2}, \sqrt{\eta_1\eta_3}, \sqrt{\eta_1}) \\
&= k_{ \{ p_1, p_2 \} } k_{ \{ p_3, p_2 \} }(\sqrt { \theta'_1 } ). \qedhere
\end{align*}
\end{proof}$\;$
\quad {\bf3.2. The 4-th multiple residue symbol.}  Let $p_1, p_2, p_3$ and $p_4$ be four prime numbers satisfying 

$$ \left\{ \begin{array}{l}
p_1 \equiv 5 \;({\rm mod} \; 8), \, p_i \equiv 1 \; ({\rm mod} \; 4)\; (i = 2, 3, 4), \\
\displaystyle{\left( \frac{p_i}{p_j} \right) = 1} \; (1 \le i \neq j \le 4), \; \mbox{$[p_i, p_j, p_k]$} = 1 \; \mbox{($i, j, k: \;$ {\rm distinct})},
\end{array} \right. 
\leqno{(3.2.1)}
$$
and we assume that the class number of $k_1 = \mathbb{Q}(\sqrt{p_1})$ is 1.\\

Let  $K$ be the field defined in Definition 3.1.7.\\\\
\quad {\bf Definition 3.2.2.}  We define the 4-{\it th multiple residue symbol} $[p_1,p_2,p_3,p_4]$ by 
$$
[p_1, p_2, p_3, p_4] = \left\{
\begin{array}{rl}
1 & {\rm if}\; p_4 \  \mbox{is completely decomposed in}\  K/\mathbb{Q}, \\
-1 &\quad \mbox{otherwise}.
\end{array}
\right.
$$
\\
We let 
$$
L := \left\{
\begin{array}{rl}
\mathbb{Q}(\sqrt{p_1}, \sqrt{p_2}, \sqrt{p_3}, \sqrt{\theta_1 \theta_2}, \sqrt{\theta_1 \theta_3})  & \mbox{if} \quad (Z,2) =1 \\
\mathbb{Q}(\sqrt{p_1}, \sqrt{p_2}, \sqrt{p_3}, \sqrt{\eta_1 \eta_2}, \sqrt{\eta_1 \eta_3}) 
& \mbox{if} \quad (Y,2) = 1 
\end{array}
\right.
$$ \\
Case $(Z,2) = 1$:  Let $\tau_1, \tau_2 , \tau_3 \in \mbox{Gal}(K/\mathbb{Q})$ be as in the proof of Theorem 3.1.8 and we let 
$$
\xi_1 := \sqrt{\theta_1 \theta_2} + \sqrt{\theta_3 \theta_4},\;
\xi_2 := \sqrt{\theta_1 \theta_3} + \sqrt{\theta_2 \theta_4}, \;
\xi_3 := \sqrt{\theta_1} + \sqrt{\theta_2} + \sqrt{\theta_3} + \sqrt{\theta_4}.
$$
Then, the subfields of $K/\mathbb{Q}$ which corresonds by Galois theory to the subgroups generated by $\tau_1, \tau_2$, $\tau_3$ and $(\tau_1 \tau_2 \tau_3 \tau_2)^2$ are $\mathbb{Q}(\sqrt{p_2}, \sqrt{p_3}, \xi_1, \sqrt{\theta_1 \theta_3}, \xi_3)$, 
$\mathbb{Q}(\sqrt{p_1}, \sqrt{p_3}, \sqrt{\theta_2}, \sqrt{\theta_3}, \sqrt{\theta_4})$, $\mathbb{Q}(\sqrt{p_1}, \sqrt{p_2}, \sqrt{\theta_1 \theta_2}, \xi_2, \xi_3)$ and $F$, respectively. By the assumption (3.2.1), $p_4$ is completely decomposed in the extension $F/\mathbb{Q}$. \\
Case $(Y,2) = 1$: We let $\tau_1, \tau_2 , \tau_3 \in \mbox{Gal}(K/\mathbb{Q})$ and $\xi_1, \xi_2, \xi_3$ be defined by replacing $\theta_i$ in the case $(Z,2) = 1$ with $\eta_i$ $(1 \leq i \leq 4)$. Then, as in the case $(Z,2) = 1$ the subfields of $K/\mathbb{Q}$ which corresponds by Galois theory to the subgroups generated by $\tau_1, \tau_2$, $\tau_3$ and $(\tau_1 \tau_2 \tau_3 \tau_2)^2$ are $\mathbb{Q}(\sqrt{p_2}, \sqrt{p_3}, \xi_1, \sqrt{\eta_1 \eta_3}, \xi_3)$, 
$\mathbb{Q}(\sqrt{p_1}, \sqrt{p_3}, \sqrt{\eta_2}, \sqrt{\eta_3}, \sqrt{\eta_4})$, $\mathbb{Q}(\sqrt{p_1}, \sqrt{p_2}, \sqrt{\eta_1 \eta_2}, \xi_2, \xi_3)$ and $F$, respectively. By the assumption (3.2.1), $p_4$ is completely decomposed in the extension $F/\mathbb{Q}$. 

Let $\frak{P}_4$ be a prime ideal in $F$ lying over $p_4$ and let $\sigma_{\frak{P}_4} = \left( \frac{K/F}{\frak{P}_4} \right) \in \mbox{Gal}(K/F)$ be the Frobenius automorphism of $\frak{P}_4$.  Note that $\frak{P}_4$ is decomposed in $K/F$ if and only if $p_4$ is completely decomposed in $K/\mathbb{Q}$. So we have, by Definition 3.2.2, 
$$
[p_1, p_2, p_3,p_4] = \left\{
\begin{array}{rl}
1 &\quad \sigma_{\frak{P}_4} = \mbox{id}_{K}, \\
-1 &\quad \sigma_{\frak{P}_4} \neq \mbox{id}_{K}. 
\end{array}
\right.
\leqno{(3.2.3)}
$$

Let $S := \{ p_1, p_2, p_3, p_4 \}$. Then, by Theorem 2.2.2, we have
\begin{align*}
G_S(2) &= \mbox{Gal}(\mathbb{Q}_S(2)/\mathbb{Q}) \\
&= \langle x_1, x_2, x_3, x_4 \mid x_{1}^{p_1 - 1}[x_1,y_1] = \cdots = x_{4}^{p_4 - 1}[x_4,y_4] = 1 \rangle.
\end{align*}
Let $\hat{F}$ be the free pro-2 group on $x_1, x_2, x_3, x_4$ and let $\pi : \hat{F}(2)  \rightarrow G_S(2) $ be the natural homomorphism. 
Since $K \subset \mathbb{Q}_S(2)$ by Theorem 3.1.10, we have the natural homomorphism $\psi : G_S(2)  \rightarrow \mbox{Gal}(K/\mathbb{Q}) $. Let $\varphi := \pi \circ \psi  : \hat{F}  \rightarrow \mbox{Gal}(K/\mathbb{Q})$. We then see that 
\[
\varphi(x_1) = \tau_1, \quad \varphi(x_2) = \tau_2, \quad \varphi(x_3) = \tau_3, \quad \varphi(x_4) = 1. 
\]
Therefore the relations among $\tau_1, \tau_2$ and $\tau_3$ are equivalent to the following relations:
\[ \begin{array}{c} 
\varphi(x_1)^2 = \varphi(x_2)^2 = \varphi(x_3)^2 = \varphi(x_1 x_3)^2 = 1,\quad \varphi(x_4) = 1,\\
\varphi(x_1 x_2)^4 = \varphi(x_2 x_3)^4 = \varphi(x_1 x_2 x_3)^4 = \varphi((x_1 x_2 x_3 x_2)^2 x_3)^2 = 1.
\end{array}
\leqno{(3.2.4)}
\]
On the other hand, by the assumption (3.2.1), we have $\overline{\mu}_2(1234) = \mu_2(1234)$.
\\
\\
\quad {\bf Theorem 3.2.5.} {\it We have}
$$
\mbox{$[p_1, p_2, p_3, p_4]$} = (-1)^{\mu_{2}(1234)}. 
$$
\begin{proof}  By (3.2.3), we have 
$$
\varphi(y_4) = \left\{
\begin{array}{ll}
1 & \; \mbox{if} \;  [p_1, p_2, p_3, p_4] = 1, \\
(\tau_1 \tau_2 \tau_3 \tau_2)^2 = \varphi((x_1 x_2 x_3 x_2)^2) & \; \mbox{if} \;  [p_1, p_2, p_3, p_4] = -1. 
\end{array}
\right.
$$
By (3.2.4), Ker($\varphi$) is generated as a normal subgroup of $\hat{F}$ by 
$$ x_1^2, x_2^2, x_3^2, (x_1x_3)^2, x_4, (x_1x_2)^4, (x_2x_3)^4, (x_1x_2x_3)^4 \; {\rm and}\; ((x_1 x_2 x_3 x_2)^2 x_3)^2
$$
 and one has
$$ \begin{array}{l}
\quad M_2( (x_1)^2 ) = (1 + X_1)^2  = 1 + X_1^2,\\
\quad M_2( (x_2)^2 ) = (1 + X_2)^2  = 1 + X_2^2, \\
\quad M_2( (x_3)^2 ) = (1 + X_3)^2  = 1 + X_3^2,  \\
\quad M_2( (x_1 x_3)^2 ) = ( (1 + X_1) (1 + X_3) )^2 \equiv 1 \; {\rm mod}\;  {\rm deg}\;  \ge 2\\
\quad M_2( (x_1 x_2)^4 ) = ( (1 + X_1) (1 + X_2) )^4   \equiv 1 \; {\rm mod}\;  {\rm deg}\;  \ge 4, \\ 
\quad M_2( (x_2 x_3)^4 ) = ( (1 + X_2) (1 + X_3) )^4 \equiv 1 \; {\rm mod} \; {\rm deg} \; \ge 4,\\ 
\quad M_2( (x_1 x_2 x_3)^4 ) = ( (1 + X_1) (1 + X_2) (1 + X_3) )^4 \equiv 1 \; {\rm mod} \; {\rm deg} \; \ge 4, \\ 
\quad M_2( ( (x_1 x_2 x_3 x_2)^2 x_3)^2 )  \\ 
\; \; \equiv 1 + X_3^2 + X_1^2 X_3 + X_1 X_3^2 + X_1 X_3^2 + X_3 X_1^2 + X_3^2 X_1
\; {\rm mod} \; {\rm deg} \; \ge 4.\end{array}
$$
Therefore $\mu_2( (1) ; *),  \mu_2( (2) ; *), \mu_2( (3) ; *), \mu_2( (12) ; *), \mu_2( (23) ; *), \mu_2( (123) ; *)$ take their values 0 on Ker($\varphi$). 
If $\varphi(y_4) = 1$, $\mu_2(1234) = \mu_2( (123) ; y_4) = 0$ by $\varphi(y_4) \in \mbox{Ker}(\varphi)$.
If $\varphi(y_4) = (\tau_1 \tau_2 \tau_3 \tau_2)^2 = \varphi( (x_1 x_2 x_3 x_2)^2 )$, we can write $y_4 = (x_1 x_2 x_3 x_2)^2 R$, where $R \in \mbox{Ker}(\varphi)$. Then comparing the coefficients of $X_1X_2X_3$ in the equality $M_2(y_4) = M_2((x_1 x_2 x_3 x_2)^2)M_2(R)$, we have
$$ \begin{array}{ll}
\mu_2(1234) & = \mu_2( (123) ; y_4) \\
            & = \mu_2( (123) ; (x_1 x_2 x_3 x_2)^2) + \mu_2( (12) ; (x_1 x_2 x_3 x_2)^2) \mu_2( (3) ; R) \\
            & + \mu_2( (1) ; (x_1 x_2 x_3 x_2)^2) \mu_2( (23) ; R) + \mu_2( (123) ; R) \\
 & = 1. 
\end{array}
$$
This yields our assertion.
\end{proof} $\;$
\\
\quad {\bf Example 3.2.6.} Let $(p_1, p_2, p_3, p_4) := (5, 8081, 101, 449)$. Then we have 
\[ \left\{ \begin{array}{l}
\theta_1 = 25 + 2 \sqrt{5} + 2 \sqrt{101},\\
\theta_2 = 25 + 2 \sqrt{5} - 2 \sqrt{101},\\
\theta_3 = 25 - 2 \sqrt{5} + 2 \sqrt{101},\\
\theta_4 = 25 - 2 \sqrt{5} - 2 \sqrt{101},
\end{array}\right.
\; 
\left\{ \begin{array}{l}
k_{\{p_1,p_2\}} = \mathbb{Q}( \sqrt{5}, \sqrt{8081}, \sqrt{241 + 100 \sqrt{5}}), \\
k_{\{p_3,p_2\}} = \mathbb{Q}( \sqrt{8081}, \sqrt{101}, \sqrt{1009+100\sqrt{101}}),\end{array}\right.
\]
and 
$$K=k_{\{p_1,p_2\}}\cdot k_{\{p_3,p_2\}}(\sqrt{25 + 2 \sqrt{5} + 2 \sqrt{101}}). $$
Then we have
$$ \left\{ \begin{array}{l}

\left( \frac{p_i}{p_j} \right) = 1\;\; (1 \leq i\neq j \leq 4),\;\;  [p_i,p_j,p_k] = 1 \;\; (i, j, k:\; \mbox{distinct}), \\
\mbox{$[p_1,p_2,p_3,p_4]$} = -1.
\end{array}
\right.
$$
In view of Example 2.1.3, this 4-tuple prime numbers may be called {\it Milnor primes}.
\begin{center}
\includegraphics[height=4cm]{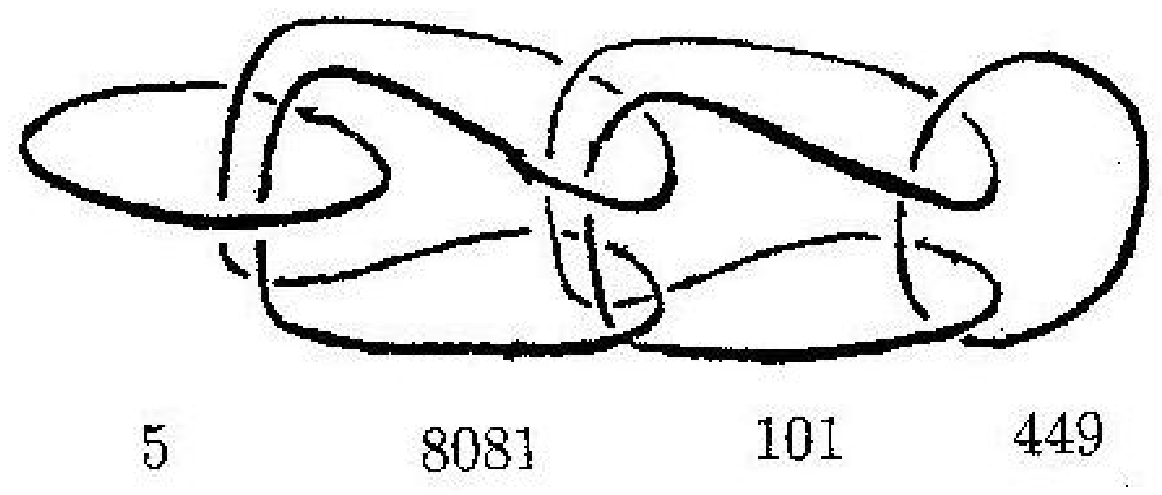}
\end{center}
Finally, two remarks are in order.\\
\\
{\bf Remark 3.2.7.}   (1) By Theorem 3.2.5, the shuffle relation for arithmetic Milnor invariants (Theorem 2.2.3 (3)) yields the following shuffle relation for the 4-th multiple residue symbol
$$ \prod_{(ijk) \in {\rm PSH}(I,J)} [p_i,p_j,p_k,p_l] = 1,$$
where $I$, $J$ are multi-indices with $|I| + |J| = 3$ and ${\rm PSH}(I,J)$ is the set of proper shuffles of $I$ and $J$, and $1\leq l \leq 4$. It is also expected that our 4-th multiple residue symbols satisfy the cyclic symmetry, although we are not able to prove it in the present paper. We hope to study the reciprocity law for the 4-th multiple residue symbol in the future.\\
\quad (2) In this paper, we are concerned only with $2$-extensions over $\mathbb{Q}$ as a generalization of R\'{e}dei's work. If a base number field $k$ contains the group of $l$-th roots of unity $\mu_l$ for an odd prime number $l$ and the maximal pro-$l$ Galois group over $k$ unramified outside a set of certain primes $S = \{ \frak{p}_1,\dots, \frak{p}_r\} \cup \{ \frak{p} | \infty \}$ is a Koch type pro-$l$ group, we can intoduce $\mu_l$-valued multiple residue symbol $[\frak{p}_1,\dots, \frak{p}_r]$ in a similar manner.
\\
\\
{\it Acknowledgements.}
I would like to thank my advisor Professor Masanori Morishita for proposing the problem studied in this paper and valuable advice.
I also thank Professor Yasushi Mizusawa for the computation of the group $N_4(\mathbb{F}_2)$ by GAP.
Finally I am grateful to the referee for useful comments.

\quad \\
\\
{\small Fumiya Amano \\
Faculty of Mathematics, Kyushu University \\
744, Motooka, Nishi-ku, Fukuoka, 819-0395, JAPAN \\
e-mail: f-amano@math.kyushu-u.ac.jp}

\end{document}